\documentclass[Afour,sageh,times]{sagej}
\usepackage{moreverb,epsfig}
\usepackage{amssymb,amsmath}
\usepackage{color}
\usepackage{amsthm}	
 \usepackage{hyperref}
 \usepackage{setspace}
\usepackage{graphicx}

\def\volumeyear{2015}

\begin{document}

\runninghead{Ji and Turkoglu}

\title{Development of a Low-Cost Experimental Quadcopter Testbed Using an Arduino Controller and Software}

\author{Ankyda Ji\affilnum{1} and Kamran Turkoglu\affilnum{2}}

\affiliation{\affilnum{1}Graduate Student, San Jose State University, Aerospace Engineering, San Jose, CA, 95192 USA\\
\affilnum{2}Assistant Professor, San Jose State University, Aerospace Engineering, San Jose, CA, 95192 USA}

\corrauth{Kamran Turkoglu, San Jose State University, Aerospace Engineering, San Jose, CA, 95192 USA}

\email{kamran.turkoglu@sjsu.edu}

\begin{abstract}
This paper explains the integration process of an autonomous quadcopter platform and the design of Arduino based novel software architecture that enables the execution of advanced control laws on low-cost off-the-shelf products based frameworks. Here, quadcopter dynamics are explored through the classical nonlinear equations of motion. Next, quadcopter is designed, built and assembled using off-the-shelf, low-cost products to carry a camera payload which is mainly utilized for any type of surveillance missions. System identification of the quadcopter dynamics is accomplished through the use of sweep data and $CIFER^{\circledR}$ to obtain the dynamic model. The unstable, non-linear quadcopter dynamics are stabilized using a generic control algorithm through the novel Arduino based software architecture. Experimental results demonstrate the validation of the integration and the novel software package running on an Arduino board to control autonomous quadcopter flights.
\end{abstract}


\maketitle

\section*{Nomenclature}
\begin{tabbing}
  XXX \= \kill
$b$ 	\> thrust factor of the propeller\\
$d$ 	\> drag factor of the propeller\\
$I_{XX}$ 	\> moment of inertia about the x-axis\\
$I_{YY}$ 	\> moment of inertia about the y-axis\\
$I_{ZZ}$ 	\> moment of inertia about the z-axis\\
$J_{TP}$ 	\> moment of inertia about the propeller axis\\
$p$ 	\> roll rate\\
$q$ 	\> pitch rate\\
$r$ 	\> yaw rate\\
$U_1$ 	\> vertical thrust factor\\
$U_2$ 	\> rolling torque factor\\
$U_3$ 	\> pitching torque factor\\
$U_4$ 	\> yawing torque factor\\
$u$ 	\> velocity in the x-axis direction\\
$v$ 	\> velocity in the y-axis direction\\
$w$ 	\> velocity in the z-axis direction\\
$\Omega$ 	\> total propellers' speed\\
$\Omega_1$ 	\> front right propeller speed\\
$\Omega_2$ 	\> rear right propeller speed\\
$\Omega_3$ 	\> rear left propeller speed\\
$\Omega_4$ 	\> front left propeller speed\\
$\phi$ 	\> roll angle\\
$\psi$ 	\> yaw angle\\
$\theta$ 	\> pitch angle\\

\end{tabbing}
\onecolumn
\section{Introduction}
Quadcopters are small rotary crafts that can be used in various environments, where they are able to maintain hover capabilities like a conventional helicopter, but are mechanically simpler and can achieve higher maneuverability. They use 4 fixed pitch propellers to control lift and a combination of propeller torques to control roll, pitch, and yaw.  Early designs had poor performance due to very high pilot workload. Current day control techniques and small sensors have increased the popularity of the quadcopter as an autonomous Unmanned Aerial Vehicle (UAV) platform.

The quadcopter is initially an unstable and underactuated plant with highly coupled and non-linear dynamics. These features make it an attractive experimental set-up and a system for controller design methodologies. There are 6 degrees of freedom to be controlled by 4 motor inputs and modeling the dynamics of the quadcopter is essential to understand its performance. For better understanding of plant behaviour, linearization of the nonlinear quadcopter model through analytical equations is essential, as investigated in \cite{bresciani2008} and \cite{balas2012}. Another aspect of the problem is the system identification portion which has been studied in  \cite{miller2007} and \cite{lee2011} and they form the basis for controller design and analysis.

In existing literature, there are many valuable works conducted on quadcopter analysis, and numerous practical applications of quadcopters ranging from disaster zone surveilence to photography, and so on. In this paper, we aim to provide a genuine approach in build and design of a prototype quadcopter purely based on low-cost, off-the-shelf products and Arduino controllers. Another unique aspect of this study is we provide a modified architecture of Arduino Mega software code, and through several modifications, we make it possible to implement many more advanced control methodologies on Arduino Mega controllers (such as adaptive, robust, optimal, sliding mode and many more). Even though existing PID based controllers on Arduino Mega boards work just fine, this removes the restriction of the code, and makes it possible to use this prototype test-bed as a research platform for the demonstration of any desired control algorithm.

As it is well known from literature, due to unstable nature of quad-copter dynamics, stabilization of the quadcopter requires moving or adding stable poles in the s-plane as discussed in \cite{sa2012}. Popular control techniques include classical Proportional Integral Derivative (PID) controllers and are studied in \cite{dicesare2009} and \cite{sa2012}. Moreover, Linear Quadratic Regulator (LQR) controllers are investigated in \cite{andreas2010} and \cite{valeria2013}. Once a stable plant is obtained, a linearized quadcopter model is used to design desired control (gains) which are then applied to the actual plant dynamics. Generally speaking, the nominal position for a quadcopter is in hover mode, where the deviation from hover is calculated using Inertial Measurement Units (IMU) and fed back into the system for control.

Modern quadcopters are heavily dependent on sensor measurements but have become more popular than ever due to the improvements in IMU and Global Positioning Systems (GPS) and their great potential for control applications. Due to heavy dependency on GPS and IMU measurements, complementary and Kalman filtering techniques are heavily used to adjust the GPS and IMU measurements to provide consistent values for the controller as shown in \cite{abas2011}. In some cases, depending on the assigned mission requirements, it is also possible to complement the sensor measurements with vision based tracking using live video during flight, where further details of such approach are provided in \cite{saakes2013,hurd2013}. However, robust controllers are required to stabilize the quadcopter from disturbances in an outdoor environment, as discussed in \cite{siebert2014}, or in fast moving references, as shown in \cite{sorensen2010} . 

In the light of these cases, this paper discusses the integration, development, build and analysis of a quadcopter platform that is aimed to operate autonomously on pre-programmed missions. The unstable, non-linear quadcopter dynamics are stabilized using a conventional (PID) controller. Included is the extraction of the linear model using analytical equations and system identification for comparison. But more importantly, with this study, we suggest a novel Arduino based architecture that enables the use and implementation of advanced level control strategies (such as adaptive control, robust multivariable control, optimal control, nonlinear control and so on) with the simple (yet effective) architecture of Arduino hardware. With this set-up, we demonstrate a fully functional and operational experimental research platform that is capable of executing autonomous guidannce and control strategies. 

The paper is organized as follows: In Section-\ref{sec:Modeling} modeling of the quadcopter dynamics is explained. Integration and the corresponding hardewares platforms are discussed in further detail in Section-\ref{sec:Integration}. Section-\ref{SSMSA} explains the state space modeling of the plant dynamics. System identification is investigated in Section-\ref{sec:SystemID} and obtained models are verified in Section-\ref{sec:validation}. With the Section-\ref{sec:conclusion} paper is concluded.

\section{Modeling} \label{sec:Modeling}
\subsection{Quadcopter Dynamics} \label{MotorDynamics}
\begin{figure}[!ht]
 \centering
 \includegraphics[width=2.5truein]{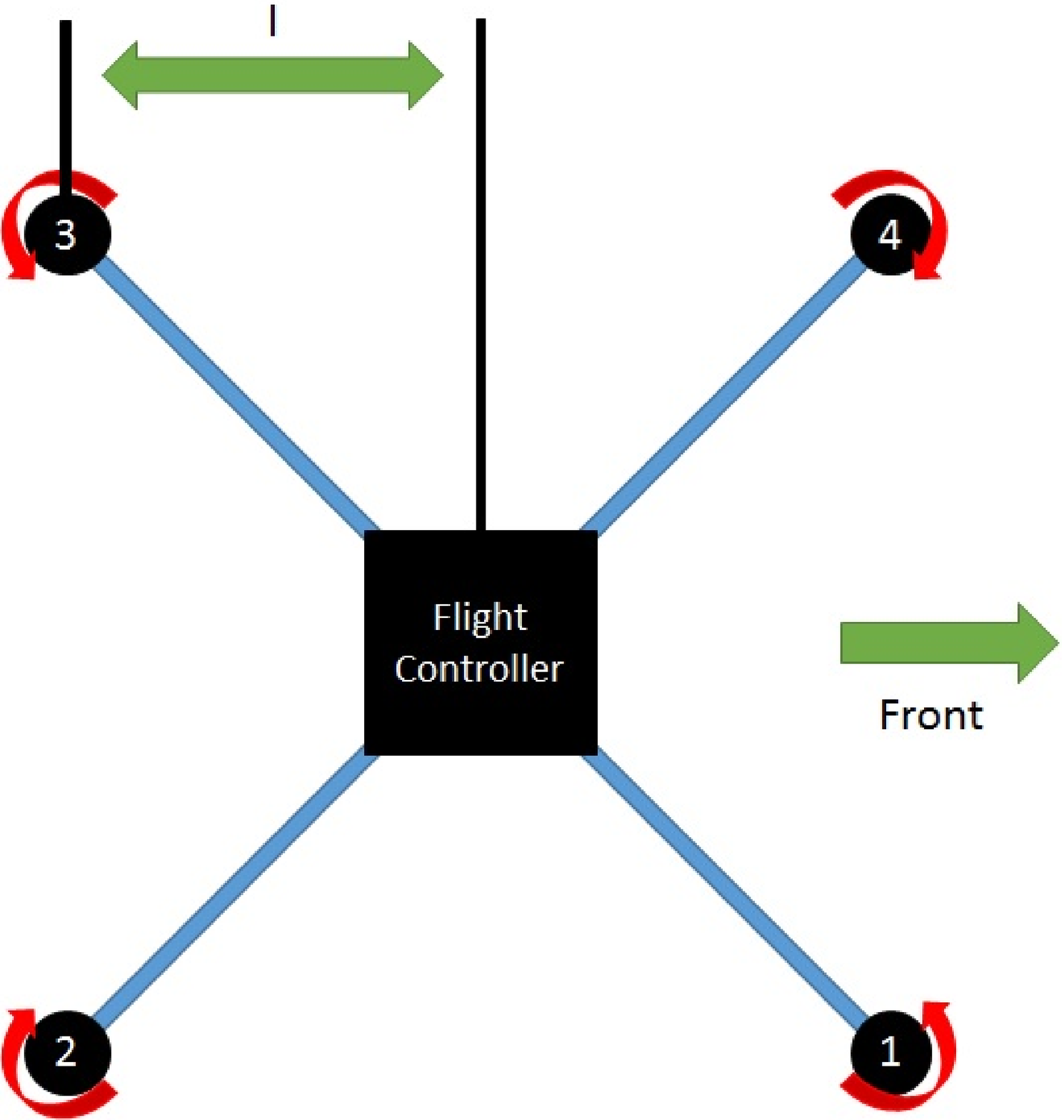}
 \includegraphics[width=3.4truein]{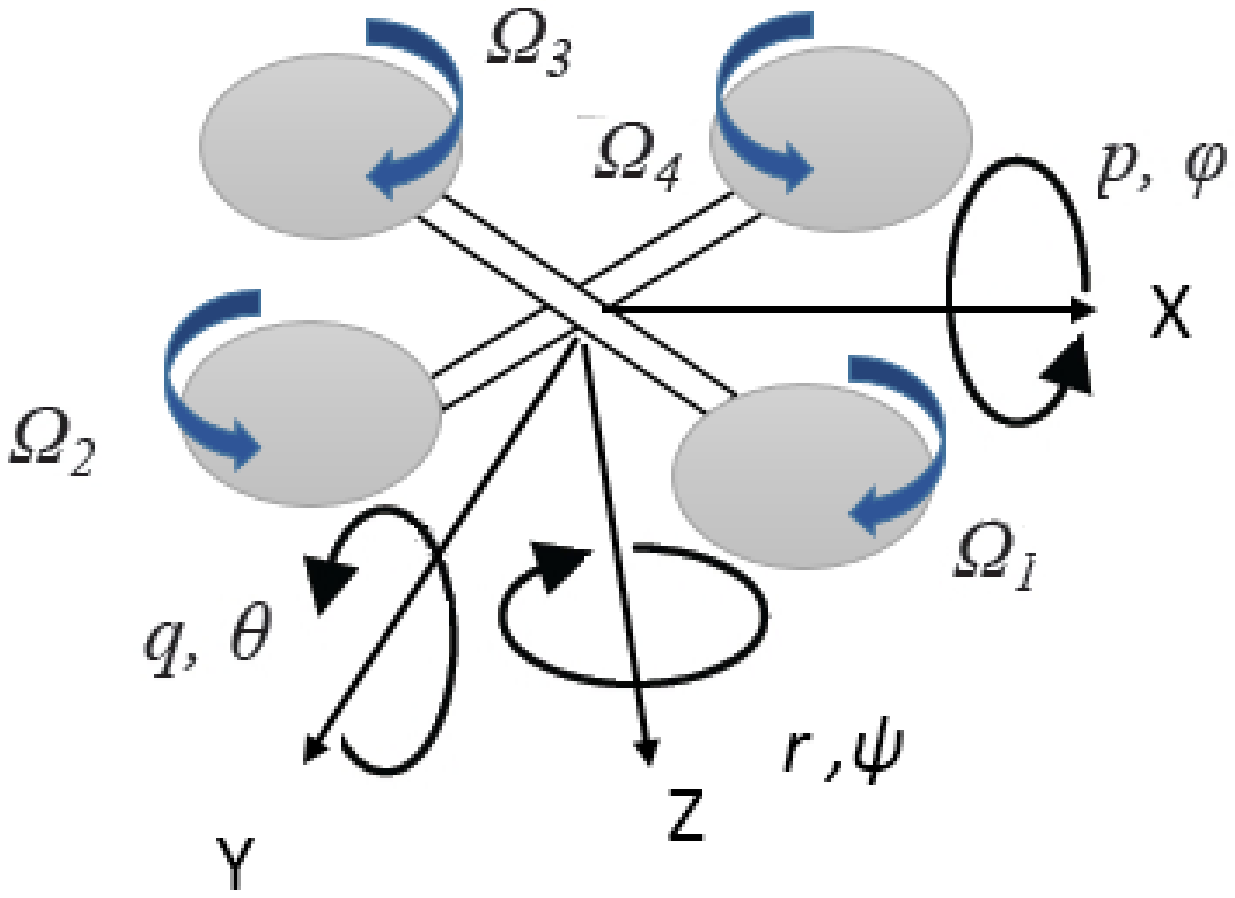}
 \caption{Quadcopter model schematic with coordinate system.}
 \label{fig:QuadDyn}
 \end{figure}

A quadcopter model consists of a cross beam structure with 4 motors on each end and collection of sets of electronic equipment. The 4 motor torques are the only inputs for a 6 degrees of freedom (DOF) system which define the quadcopter as an underactuated system. Without a controller to compensate for underactuation, there are 2 states that cannot be directly commanded. This, eventually, will cause a drift to undesired values over time. The motors are stationary and do not have any mechanical linkages to change the blade pitch. In that sense, the quadcopter utilizes a combination of the four motor torques to control all the states. The testbed is designed using the X-formation shown in diagrams of Figure \ref{fig:QuadDyn}. One pair of propellers ($\Omega_1$ and $\Omega_3$) rotate clockwise (CW) while the other pair of propellers ($\Omega_2$ and $\Omega_4$) rotate counter-clockwise (CCW).

To command throttle, all four propellers must rotate at the same speed which provides a vertical force in the z-axis. If each propeller provides a quarter of the weight in thrust, the quadcopter will hover. Rolling motion is generated by either increasing or decreasing the torque in pair of motors on the left side ($\Omega_3$ and $\Omega_4$) while applying an opposite increase or decrease to the right pair of motors ($\Omega_1$ and $\Omega_2$). This produces a torque in the x-axis which creates the rolling motion. Pitching motion is generated by increasing/decreasing the front motors ($\Omega_1$ and $\Omega_4$) while applying the opposite action (increase/decrease) to the rear motors ($\Omega_2$ and $\Omega_3$). This combination produces a torque in the y-axis which creates a pitching motion. Yawing motion is generated by increasing/decreasing the CW motor pair while applying the opposite action  (increase/decrease) to the CCW motor pair. This combination produces a torque in the z-axis which creates a yawing motion.

\subsection{Equations of Motion} \label{EquationsofMotion}
The quadcopters' non-linear, coupled equations of motion (EoMs) have been analyzed extensively in literature and are summarized below for convenience, where further detailed derivations are provided in \cite{bresciani2008} and \cite{balas2012}. These EoMs are derived by applying Newton's $2^{nd}$ Law to the quadcopter body. Some of the basic assumptions include that i) the quadcopter is a rigid body and ii) it is symmetrical along the x and y axes.
 
\begin{equation}\label{eq:eoms}
\begin{split}
&\dot{u} = (v r - w q) + g s_\theta \\
&\dot{w} = (w p - u r) - g c_\theta s_\phi \\
&\dot{w} = (u q - v p) - g c_\theta s_\phi + \frac{U_1}{m} \\
&\dot{p} = \frac{I_{YY} - I_{ZZ}}{I_{ZZ}} q r - \frac{J_{TP}}{I_{XX}} q \Omega + 
\frac{U_2}{I_{XX}} \\ 
&\dot{q} = \frac{I_{ZZ} - I_{XX}}{I_{YY}} p r - \frac{J_{TP}}{I_{YY}} p \Omega + 
\frac{U_3}{I_{YY}} \\
&\dot{p} = \frac{I_{XX} - I_{YY}}{I_{ZZ}} p q - \frac{J_{TP}}{I_{ZZ}}
\end{split}
\end{equation}

The outputs of the EoMs are translational velocities u, v, w; rotational velocities p, q, r; positions x, y, z; and the attitude angles $\phi$, $\theta$, and $\psi$. These outputs are calculated by integrating the EoMs, given in Eq.\eqref{eq:eoms}. The inputs of the EoMs are the propeller speed inputs where $U_1, U_2, U_3, U_4$ are associated with throttle, roll, pitch and yaw respectively. Here, $\Omega$ is the sum of the propellers rotational speed. These inputs are functions of the propellers rotational speed $\Omega_{1-4}$, where lift and drag factors of the propeller blade (b and d respectively) and length l. Here,  lift and drag factors of the propeller blade (b and d respectively) are calculated from the Blade Element Theory \cite{bresciani2008} . From the quadcopter dynamics discussed in 
Section \ref{MotorDynamics}, the inputs can be expressed as shown in Eq.(\ref{MotorEqns}).

\begin{equation}
\begin{split}
&U_1 = b (\Omega_1 ^2 + \Omega_2 ^2 + \Omega_3 ^2 + \Omega_4 ^2) \\
&U_2 = l b (-\Omega_1 ^2 - \Omega_2 ^2 + \Omega_3 ^2 + \Omega_4 ^2) \\
&U_3 = l b (\Omega_1 ^2 - \Omega_2 ^2 - \Omega_3 ^2 + \Omega_4 ^2) \\
&U_4 = d (-\Omega_1 ^2 + \Omega_2 ^2 - \Omega_3 ^2 + \Omega_4 ^2) \\
&\Omega = b (\Omega_1 - \Omega_2 + \Omega_3 - \Omega_4)
\label{MotorEqns}
\end{split}
\end{equation}

\subsection{Linearization of Non-Linear EoMs} \label{SimMod}
A linearized model of the EoMs are desired to analyze quadcopter behavior at an operating point where all states are effectively zero. For the quadcopter, the operating point is selected as the hover condition where the motors provide enough thrust to counteract the force of gravity. Matlab is used to linearize the model dynamics,

\begin{figure}[!ht]
 \centering
 \includegraphics[width=4truein]{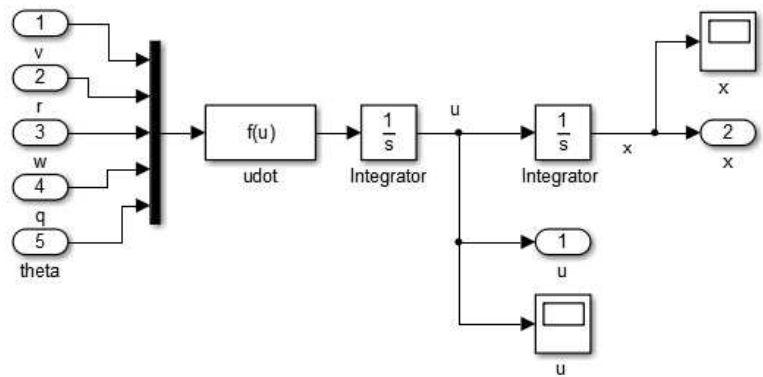}
 \caption{Simulink model of $\dot{u}$ equation, where the remaining equations follows the same.}
 \label{fig:uDotSimulink}
 \end{figure}

where Figure \ref{fig:uDotSimulink} demonstrates one of the 6 non-linear EoMs that are modeled in Simulink. The function block contains the $\dot{u}$ equation. This block is integrated twice to obtain the velocity u and position x. The control system toolbox embedded linmod function is utilized to extract the linear state space model (SSM). The same procedure is repeated for the remaining equations.

\section{Quadcopter Platform Hardware Components} \label{sec:Integration}
\begin{figure}[!ht]
 \centering
 \includegraphics[width=3truein]{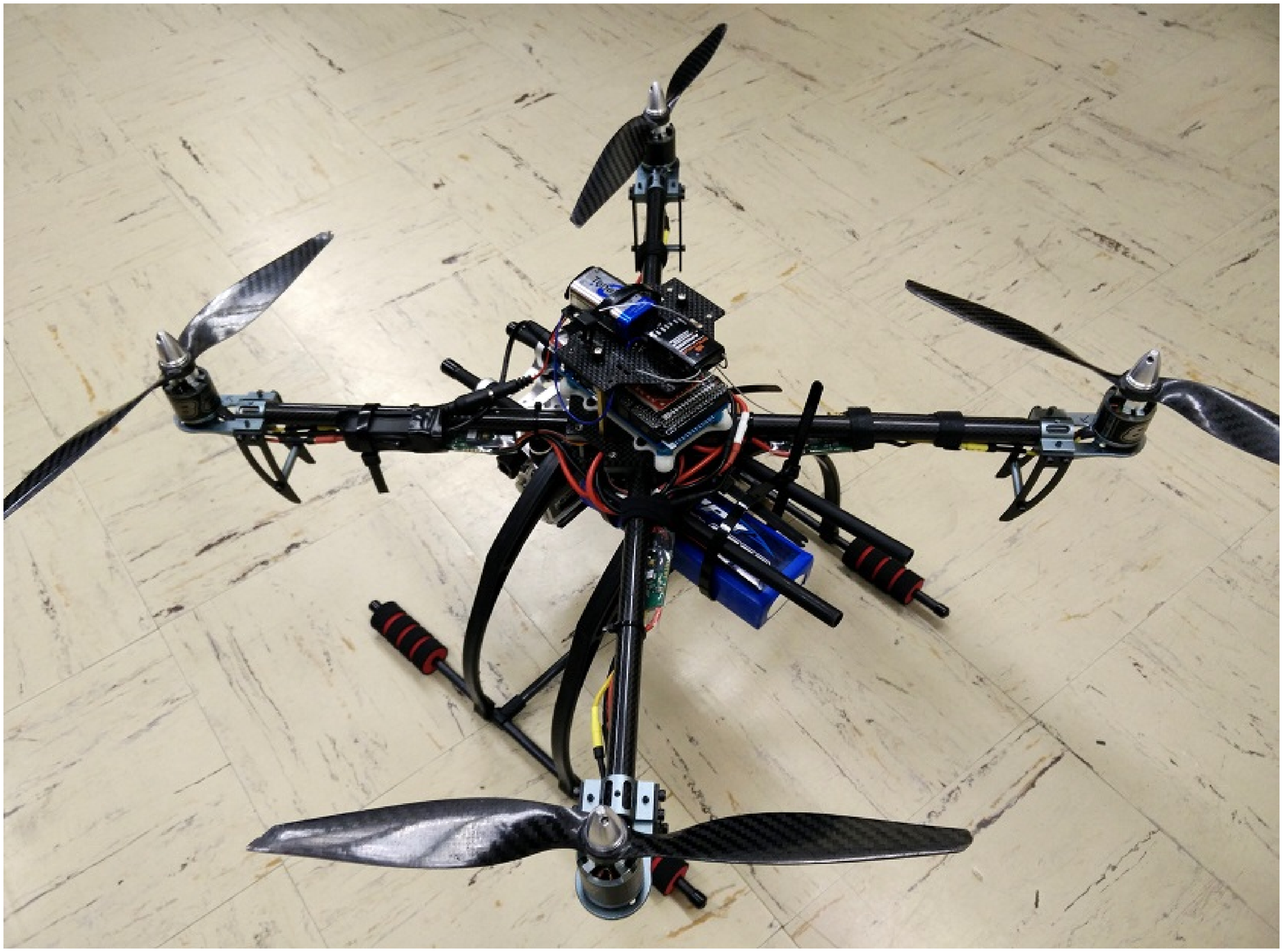}
 \includegraphics[width=3truein]{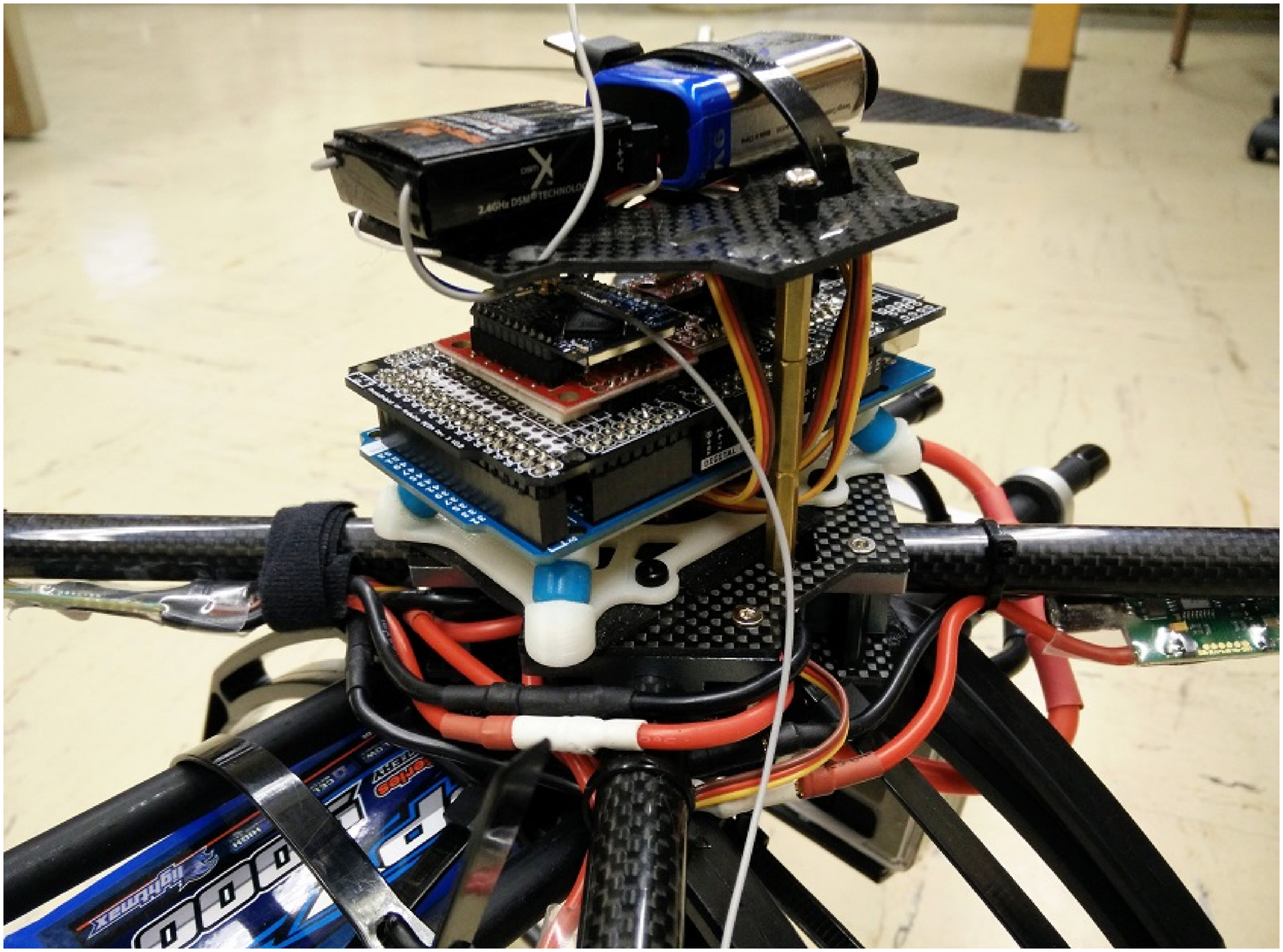}
 \caption{Assembled quadcopter testbed.}
 \label{fig:QuadPhoto}
 \end{figure}

The quadcopter testbed is built and assembled from scratch. A standard thrust to weight ratio of 2 and above is kept for maneuverability of the quadcopter. This directly results in a desire for small, lightweight components. Utilization of low cost components is one important figure of merit but the reliability of each part is also taken into account extensively. The total cost of the quadcopter is around \$388 (without a GoPro camera) and a breakdown of the off-the-shelf components can be seen in Table \ref{table.cost}. The testbed (in its current condition) can be seen in Figure \ref{fig:QuadPhoto}. 

\begin{table*}[!ht]
\centering
\caption{Quadcopter components and cost.}
\begin{tabular}{|l|l|l|}
  \hline
  Item List			& Specifics				& Price (\$) \\ \hline
  Frame				& Turnigy Talon V2.0			& 54 \\
  Arduino			& Arduino Mega 2560 R3			& 38 \\
  IMU				& Geeetech 10DOF			& 10 \\
  GPS				& GlobalSat EM 406a			& 38 \\
  ESC				& Afro 30A				& 13 \\
  Motor(4 motors)		& NTM 28-30S 800kv 300W short shaft	& 15 \\
  Propeller(2 pairs)		& 12x6 carbon fiber			& 9 \\
  Xbee				& Xbee Pro 900 RPSMA			& 20 \\
  Power Distribution Board	& Hobbyking breakout cable		& 4 \\
  6 Channel Receiver		& Spektrum AR6255			& 49 \\
  Landing Gear			& Amazon				& 24 \\
  Battery			& 3S1P and 4S1P 5000mAh 30C		& 21 \\
  \hline
                  ~~            & {\bf Total} & 295 \\
  \hline
\end{tabular}
\label{table.cost}
\end{table*}

\subsection{Motors and Speed Controllers}
\begin{figure}[!ht]
 \centering
 \includegraphics[width=2truein]{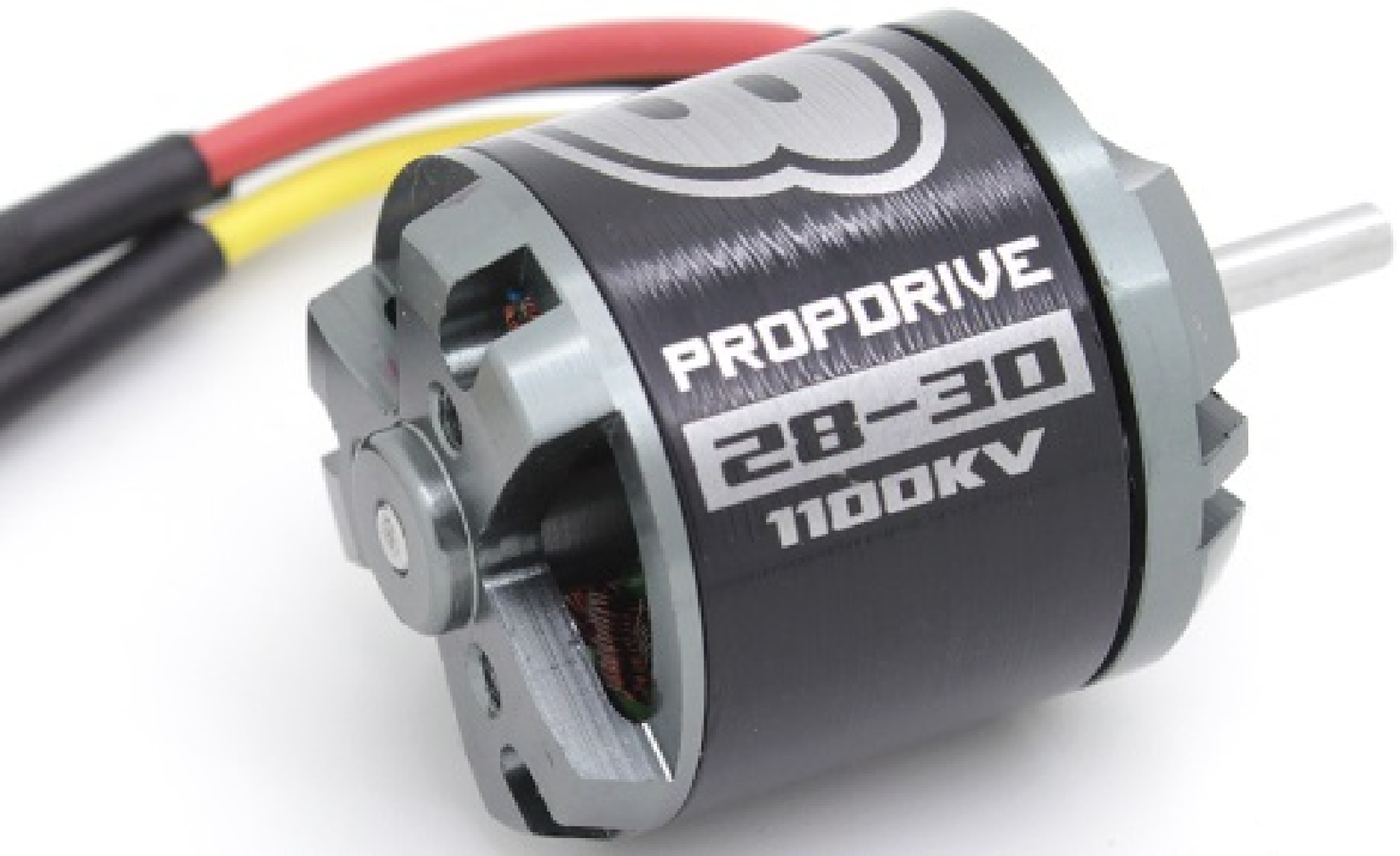}
 \includegraphics[width=2truein]{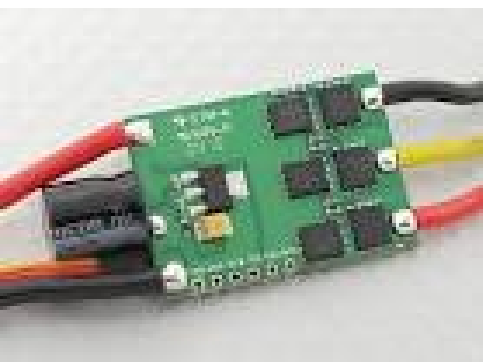}
 \caption{NTM BLDC Motor and Afro ESC}
 \label{fig:MotorESC}
 \end{figure}

To keep the thrust to weight ratio above 2, each motor and propeller combination must provide at least one half of the total weight in thrust. In order to accomplish this, NTM Prop Drive 28-30 800kv BLDC motors (seen on the left in Figure \ref{fig:MotorESC}) were selected to provide thrust using 12x6 carbon fiber propellers. The low kv rating provides high torque at lower rotational speeds to push the large props. Initial testing of these motors, using a custom thrust stand, show that each motor and propeller combination provide 900 grams of thrust at around 150 Watts. Further thrust testing at maximum output at 300 Watts show an increase in thrust to 1.3 kilograms of thrust. This results in a thrust to weight ratio of 2.8 which leaves extra thrust capacity for a camera or other accesories to be added on in the future. Afro electronic speed controllers (ESCs) are used to control the motor and propeller rotational speed and can be seen on the right in Figure \ref{fig:MotorESC}. The stock ESCs have an update rate of 50Hz which are improved by using SimonK firmware to obtain higher update rates up to 400Hz to command rotational speed.

\subsection{Battery}
Two different batteries were used in combination with the motors to test the performance and flight time of the quadcopter. Originally, a 3 cell 1 series (3S1P) 5000 milliAmphour (mAh) battery with 30C discharge rate was used. The 3S1P showed that the motors required over 50\% power to achieve takeoff. This did not meet the power to weight ratio of greater than 2. The battery was upgraded to a 4S1P 30C 5000mAh battery. It is heavier by 100 grams, but increased the power to weight ratio to meet the requirement. 

\subsection{Instrumentation}
\begin{figure}[!ht]
 \centering
 \includegraphics[width=2truein]{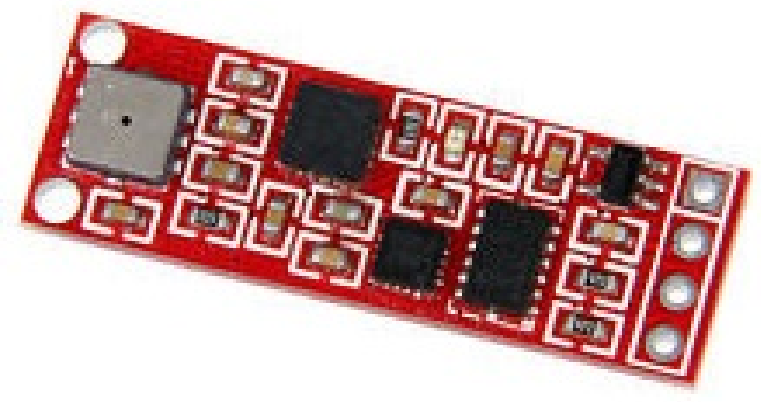}
 \includegraphics[width=2truein]{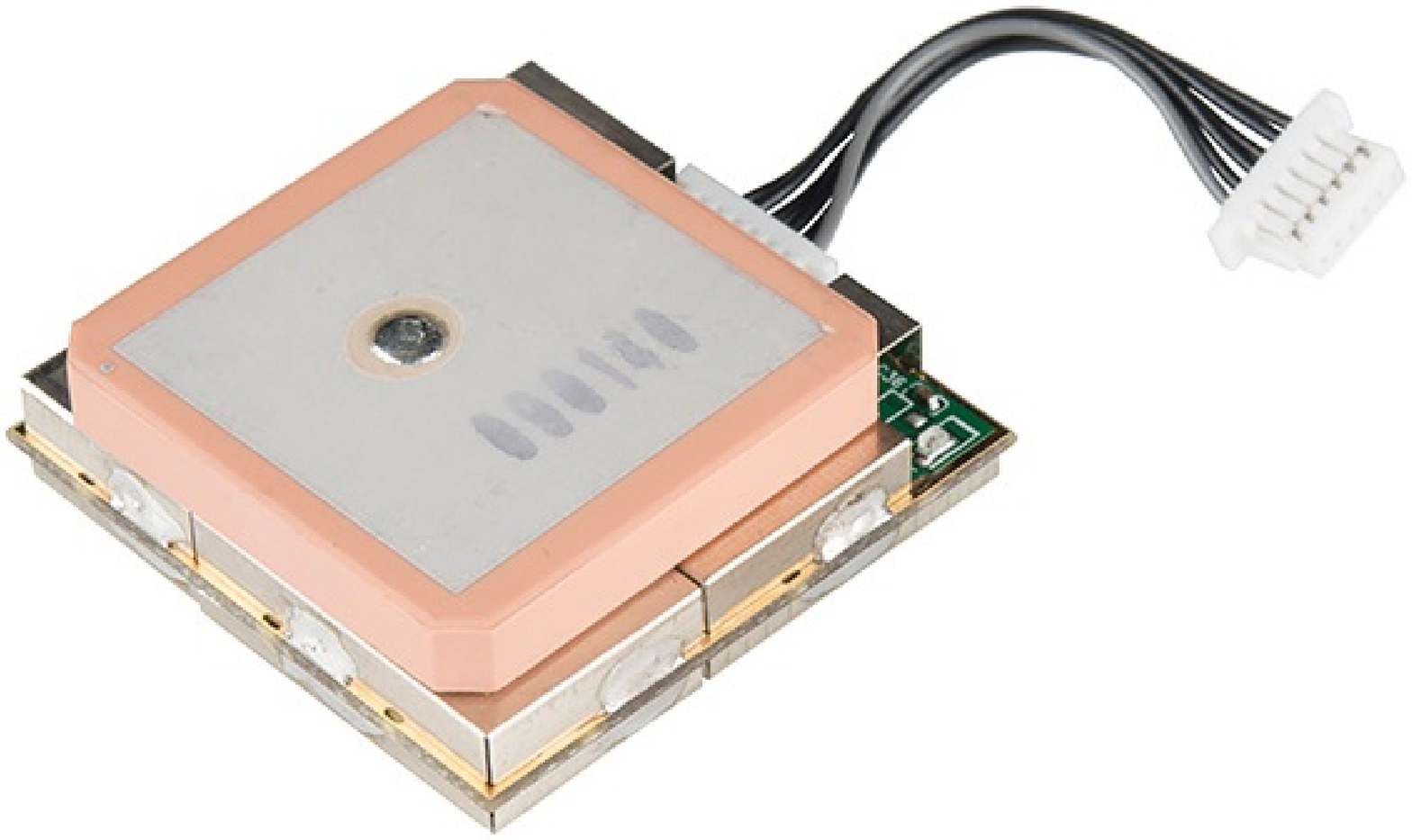}
 \caption{Geeetech 10 Degree of Freedom IMU and GlobalSat GPS.}
 \label{fig:IMUGPS}
 \end{figure}

The onboard IMU is a Geeetech 10DOF board which can be seen in on the left in Figure \ref{fig:IMUGPS}. This IMU includes a 3 axis ADXL345 accelerometer with 13-bit resolution, a 3 axis L3G4200D gyro with 16-bit resolution, a 3 axis HMC5883L magnetometer with an 12-bit resolution, and a BMP085 barometer which can achieve 0.03hPa accuracy. The GlobalSat EM-406a Global Positioning System (GPS) is included to complement the IMU measurements with positional data. It can be seen on the right in Figure \ref{fig:IMUGPS}. This GPS has an accuracy in position measurements of 5 meters with Wide Area Augmentation System (WAAS) enabled and consumes 44mA at 4.5V~6.5V.

\subsection{Receiver and RC controller}
\begin{figure}[!ht]
 \centering
 \includegraphics[width=2truein]{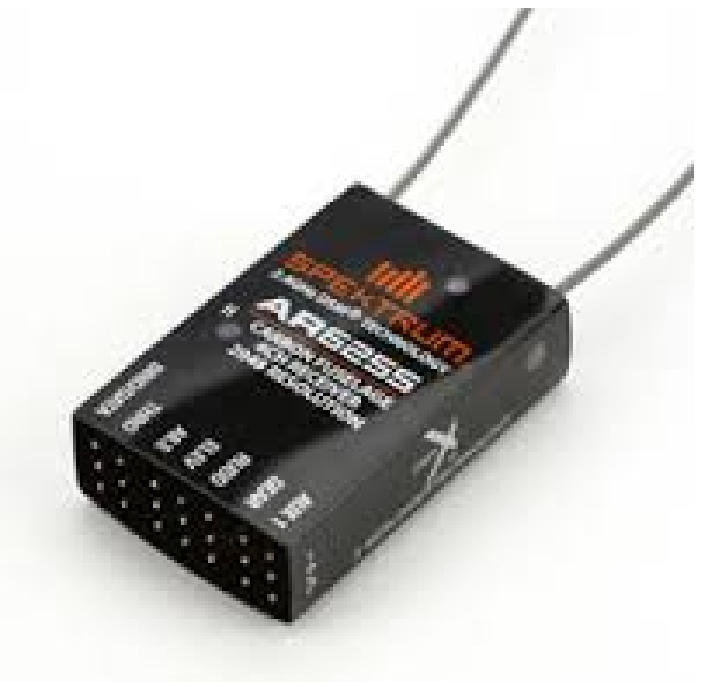}
 \includegraphics[width=2truein]{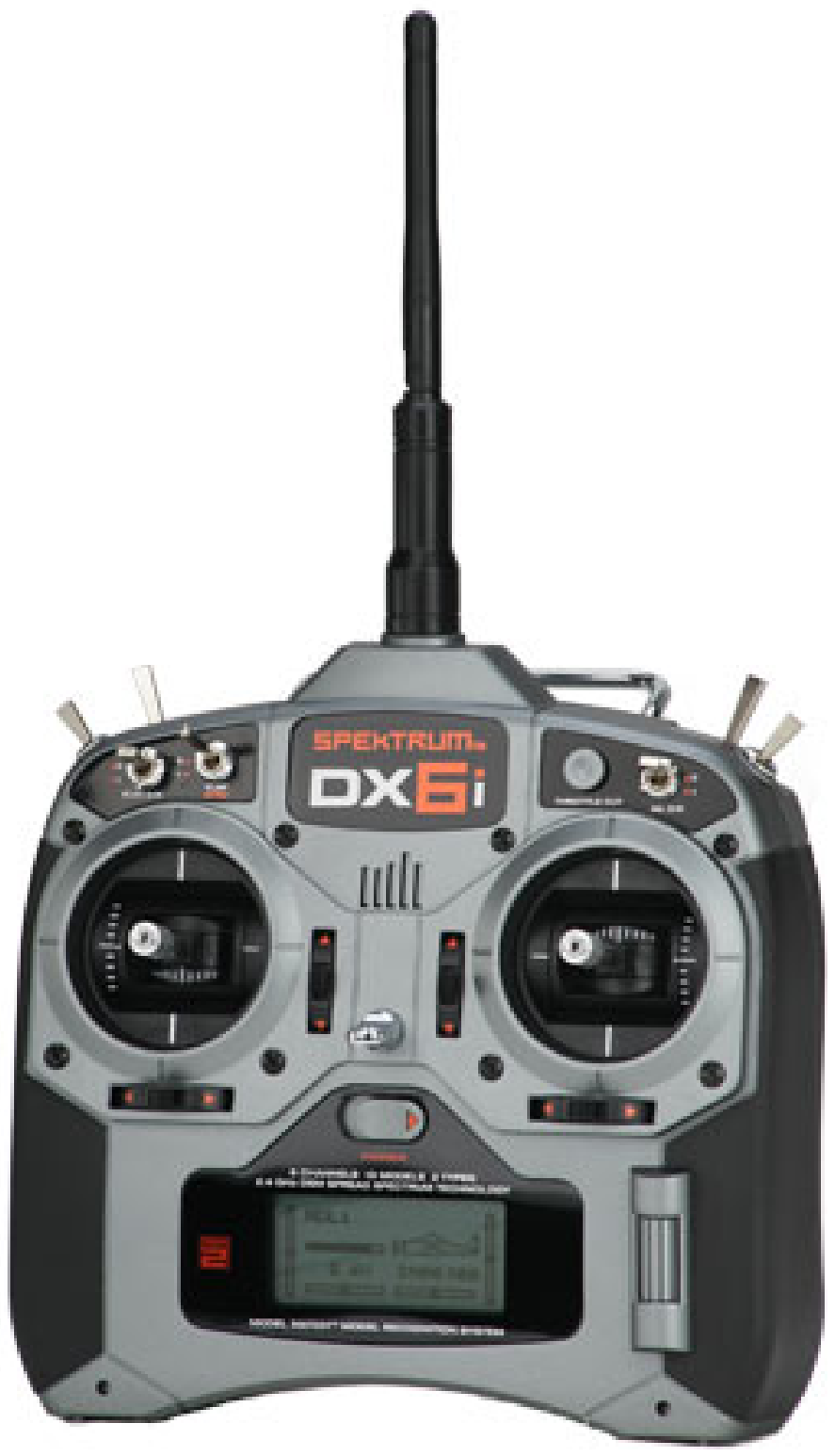}
 \caption{6 channel receiver and transmitter from Spektrum}
 \label{fig:RXDX6i}
 \end{figure}

The Spektrum AR6255 six channel receiver (seen on the left in Figure \ref{fig:RXDX6i} is used in combination with a DX6i RC controller(seen on the right in Figure \ref{fig:RXDX6i}. An input signal is transmitted from the DX6i to the paired receiver. The receiver signal is read by the microcontroller which decides to send the original signal or calculate a closed loop signal for the motors. This receiver/controller pair is used so that a pilot can input commands for initial testing and for recovery if an accident occurs. Four of the six channels are used for throttle, roll, pitch, and yaw. The 5th channel is used as a trigger switch to initiate a programmed script. The last channel is used as a kill switch which will turn off the closed loop controller for safe ground handling.

\subsection{Data Logging}
\begin{figure}[!ht]
 \centering
 \includegraphics[width=2truein]{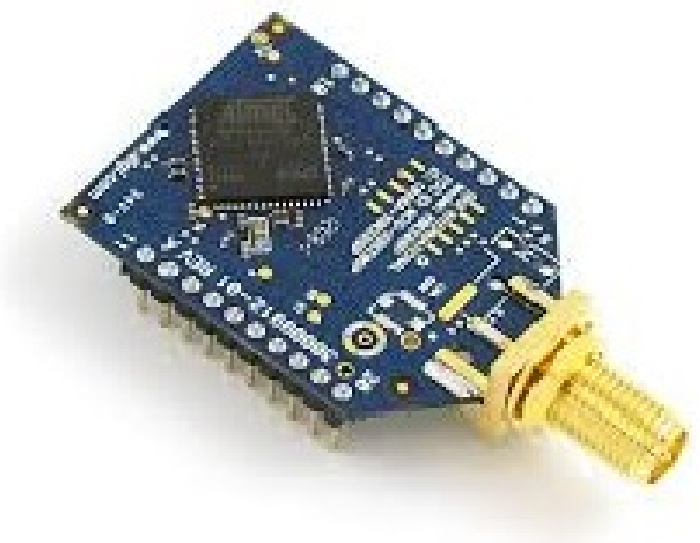}
 \includegraphics[width=2truein]{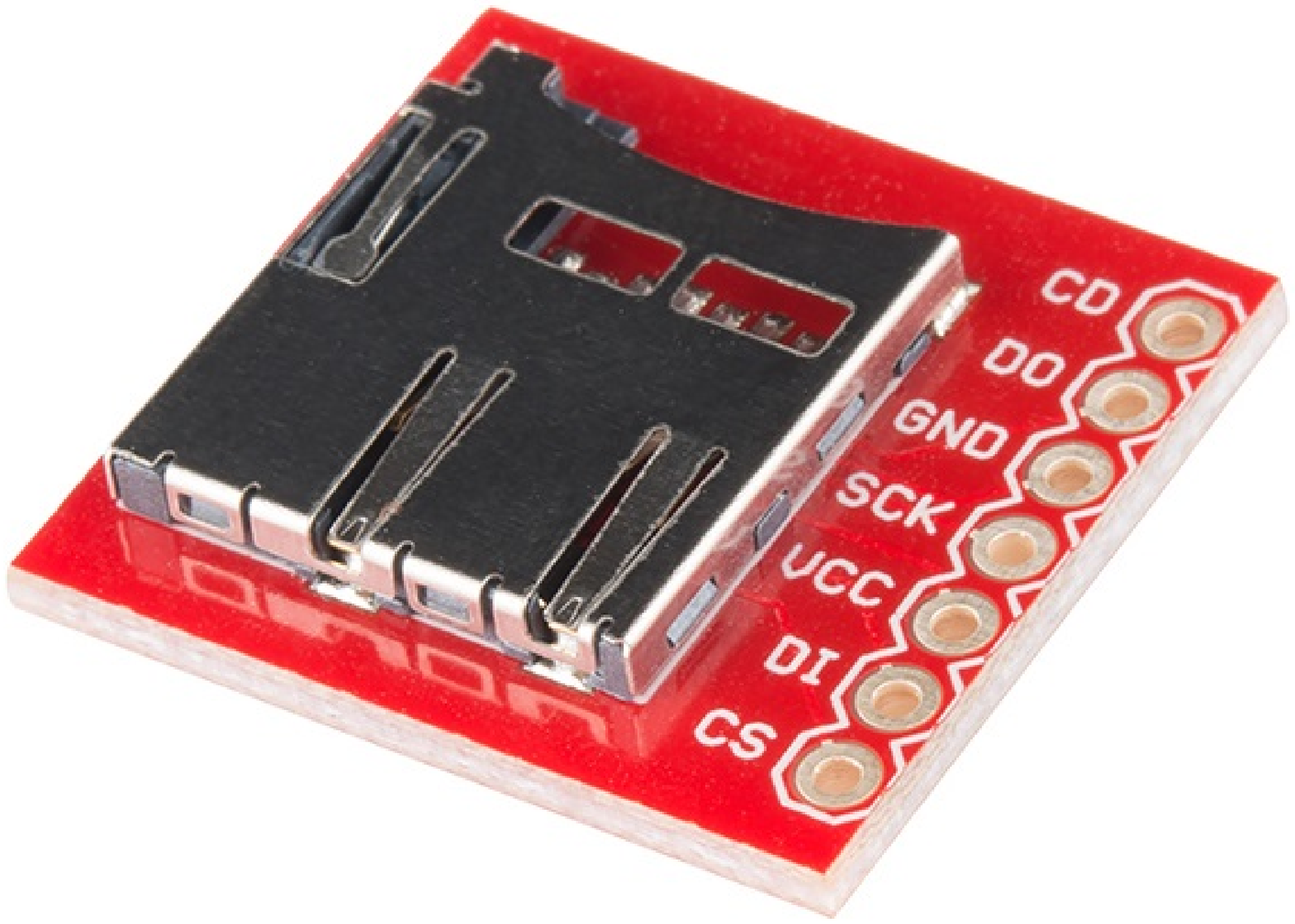}
 \caption{Xbee Pro 900 and Sparkfun MicroSD breakout board}
 \label{fig:XBeeSD}
 \end{figure}

Data from the instrumentation about the quadcopter states are kept by logging the values. Data logging was originally done using the Xbee Pro 900 RPSMA seen on the left in Figure \ref{fig:XBeeSD}. This model is more expensive but it runs on the 900MHz frequency which does not interfere with the receiver (2.4GHz) and the GPS (5.8GHz). One module is onboard the quadcopter to tramsit packets to the second module which is connected through USB to a laptop. 

The two Xbee modules are connected to each other using the XCTU program which is freeware from the company Digi. The process involves first adding the Xbee module that is connected through USB. Then the program allows you to search for other Xbee devices on the same frequency. Once the second Xbee is discovered, the program allows you to pair them to send or receive packets from one another. The Xbee option was replaced with a micro SD card due to the Xbee having low reliability in transmission from packet loss. The Xbee is currently used as a quick debug where the data can be viewed from the laptop instead of removing the micro SD card from the microcontroller and using the USB adapter.

A micro SD card breakout board from Sparkfun (seen on the right in Figure \ref{fig:XBeeSD} ) was purchased and integrated into the quadcopter. Since this SD breakout is connected directly to the microcontroller, it collects data more reliably than using the Xbee which would lose a line of data unexpectedely. The microcontroller gathers data and logs into a comma separated value (.csv) file on the micro SD card.

\subsection{Frame and Landing Gear}
The Turnigy Talon Quadcopter (V2.0) fiber frame is purchased as the frame for the Quadcopter. The Talon kit comes with carbon fiber parts, which provides higher strength and about 130 grams of weight savings over other commercially available quadcopter frames. Tall landing gear is used to provide clearance for a camera and gimbal pair which will be attached below the frame.  

\subsection{Microcontroller Setup}
\begin{figure*}[!ht]
 \centering
 \includegraphics[width=6truein]{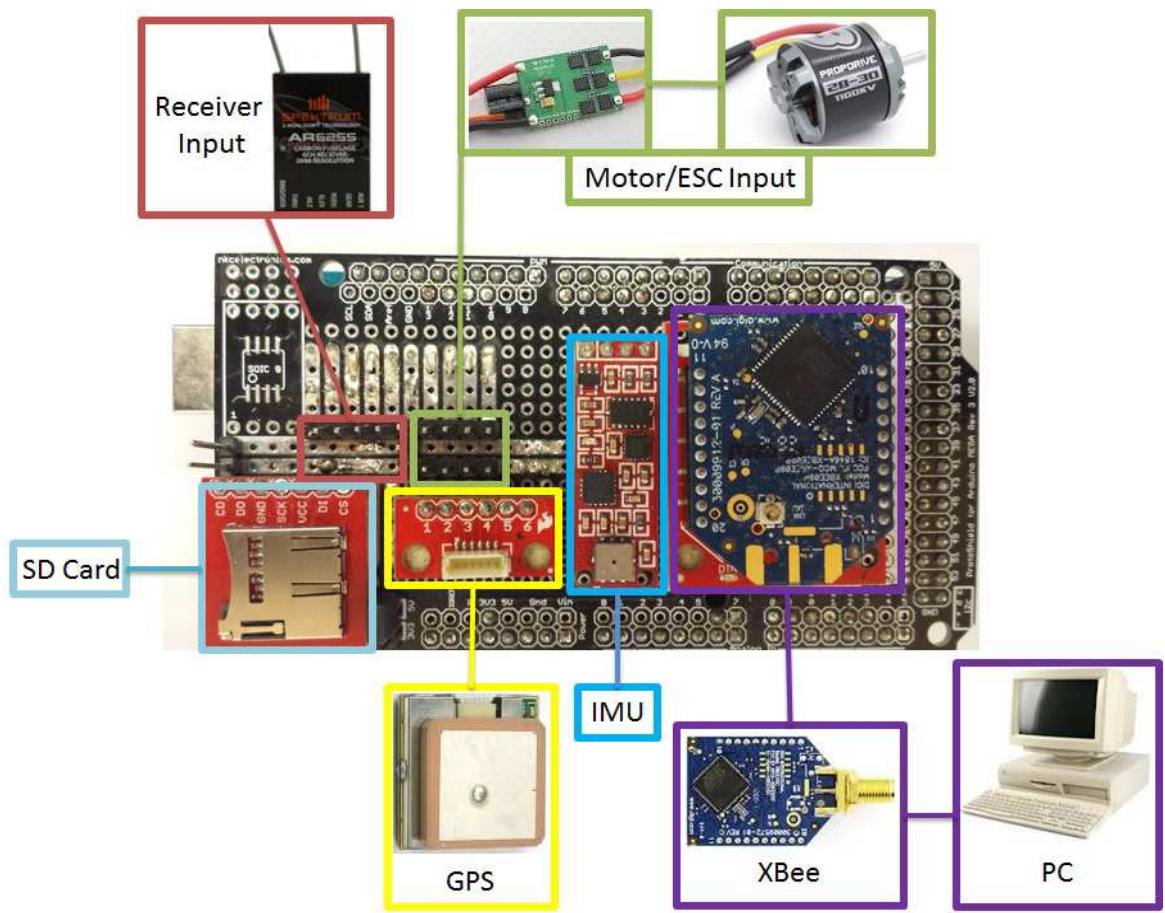}
 \caption{Custom made quadcopter shield for Arduino Mega 2560}
 \label{fig:QuadShield}
 \end{figure*}

In this study, the Arduino Mega 2560 is used as the microcontroller due to its inexpensive nature and relatively powerful characteristics. It is simple to program and its open source nature, with extensive documentation, provides tremendous benefits. The Arduino Mega operates the ATmega2560 chip microcontroller and runs at 16MHz clock speed. The Arduino Mega contains 54 digital input/output pins that can be used to control the ESCs and receiver, an I2C bus, and three pairs of TX/RX pins can be used for communication by serial transmission. Instead of connecting all the different modules to the Arduino board with scattered wiring, a custom sautered shield was designed and created using the Mega ProtoShield from AdaFruit as the base. This shield incorporates all the modules on the top and hides all the wiring underneath. A a diagram of the shield can be seen in Figure \ref{fig:QuadShield}

\begin{figure}[!ht]
 \centering 
 \includegraphics[width=3.5truein]{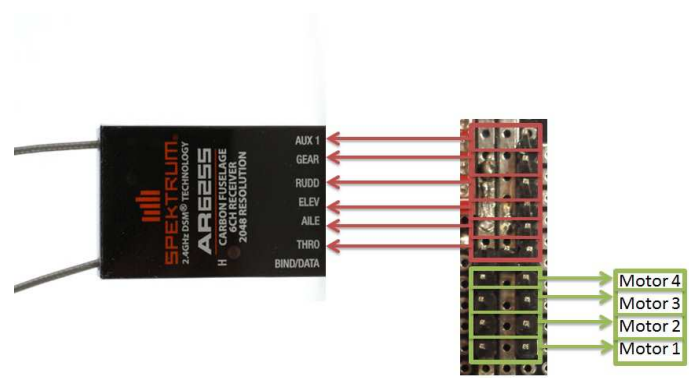}
 \caption{Detailed figure of motor and receiver of quadcopter shield}
 \label{fig:QuadShieldZoomed}
 \end{figure}

All six lines from the receiver are used and connected to the shield. The shield then connects these lines to the pins A8-A13 on the Arduino. The receiver input signals are then read using Pulse Width Modulation (PWM) by the triggering of interrupts which will be discussed further in Section \ref{PWMmapping}. The receiver input signals are then mapped into output ESC values. The output signals are sent from Arduino to the ESC which will then control the motor speed. A zoomed in diagram of the receiver and motor pins are shown in Figure \ref{fig:QuadShieldZoomed} which shows the location of each pin location on the shield.

The IMU is connected to the Arduino Mega through the I2C ports. Each module on the IMU has a unique I2C address to separate the measurement data. GPS, telemetry, and video footage will provide data to be transmitted over the serial TX/RX pins. The Arduino will process and calculate the measurement data which will be saved on the micro SD card or alternatively transmit all of it over the wireless XBee. A camera is implemented underneath which can be used to stream video wirelessly to a laptop or phone. A summary of the arduino pins used for each component is shown in the Table \ref{table.ArduinoPins}.

\begin{table}[!ht]
\centering
\caption{Arduino pin configuration.}
\begin{tabular}{|l|l|l|}
  \hline
  Component 			 & Type of Pin & Pin value(s) \\ \hline
  Electronic speed controllers   & Digital I/O & 3, 5, 6, 7 \\
  Spektrum 6 channel receiver    & Analog I/O  & A8 - A13 \\
  Inertia measurement unit 	 & I2C 	       & SCL, SDA \\
  EM406a GPS 		 	 & TX/RX       & 16, 17 \\
  XBee telemetry 		 & TX/RX       & 18, 19 \\
  Camera		 	 & TX/RX       & 14, 15 \\
  SD Card			 & SPI	       & 50 - 53 \\
  \hline
\end{tabular}
\label{table.ArduinoPins}
\end{table}

\subsection{Estimated Parameters}
Some of the quadcopter testbed parameters are listed in Table \ref{table.parameters} for convenience. The inertia values are calculated using moment of inertia equations and the parallel axis theorem which can be found in any civil engineering textbook. Values for b and d are estimated using Blade Element Analysis. These parameters are used in the EOMs to estimate an initial controller for stability which is discussed further in section \ref{SSMSA}

\begin{table}[!ht]
\centering
\caption{Calculated quadcopter platform parameters.}
\begin{tabular}{|l|l|l|}
  \hline
  Parameter & Value 		 & Units \\ \hline
  mass      & 1800 		 & grams \\
  I$_{XX}$  & 7.06*10$^{-3}$     & kgm$^2$ \\
  I$_{YY}$  & 7.06*10$^{-3}$     & kgm$^2$ \\
  I$_{ZZ}$  & 7.865*10$^{-3}$    & kgm$^2$ \\
  J$_{TP}$  & 14.2*10$^{-4}$     & kgm$^2$ \\
  b         & 4.5*10$^{-4}$      & meters \\
  d         & 1.8*10$^{-5}$      & meters \\
  l         & 8.25		 & meters \\
  \hline
\end{tabular}
\label{table.parameters}
\end{table}

\section{Software Development}\label{sec:software}
There are many different quadcopter controllers that use Arduino based software (Ardupilot, MultiWii, etc..), but one major trait in those controllers is that they are difficult to modify and lack modularity. This becomes problematic when researchers want to use such test-beds for the application of more advanced (and sophisticated) control methodologies, such as adaptive control, mu-synthesis, and so on. The reason why many of these systems  (like Ardupilot, MultiWii, etc..) are difficult to modify is that they are very limited in terms of their own, specialized libraries and specifically named functions/values. Currently, there is also not enough documentation provided to explain each software, due to its open-source nature. A specific example can be found in Ardupilot whose stabilization algorithm uses a PID library which is calculated using values that are found in a different library making it a very round about process. 

In this research, for the sake of modularity and applicability, the Arduino code was written from scratch. Emphasis was put on simplicity by providing important values like roll and pitch angles in the main script and by using better notation scheme that relates a value to what it actually means. This made it much more simple to modify and to apply different stability and control algorithms as well as being applied to hands on experience for classroom/teaching usage. 

Most arduino codes are based on using libraries to provide utility functions that can be used in the main sketch. A library is essentially used to provide the code with a complete procedure without clumping up the main sketch in our novel code modification. The fundamental calculations of our novel modified code are written in libraries and rely on different libraries to provide quadcopter state information. This makes it very simple to apply a formula in the main script just by naming the function and providing the data that is used in the calculation. In our novel code modification, IMU data is read through a custom library while the GPS makes use of the open source TinyGPS library. GPS and IMU data are all read through the library and provide measurements for the main Arduino code to calculate PID values for stability. Since the inputs for calculations are all from libraries, it is very simple to apply (if necessary) other sensor inputs to provide the same information, which makes the Arduino code modular. Currently the arduino code uses a modified version of the open source PID library by \cite{Beauregard2012}, but can be easily swapped out for other controller schemes (such as LQR, Adaptive control, Robust control ... etc.) by creating new libraries with the corresponding formulas/calculations and including them in the main script.

Some of the known limitations of the modified arduino system include sampling time, telemetry, and accurate measurement data. Sampling time is relatively difficult to keep consistent because it is a function of how many calculations must be completed. If more complicated controllers or additional modules are added to the system, more computing power is required which will lower the sampling time to a point where the quadcopter's stability is endangered. The current code runs at 100[Hz] which leaves enough safety factor for calculations, and therefore the stability. Previous hardware testing had showed unstable behavior below 20[Hz]. 

With all these in mind this code provide the flexibility and modularity of not only improved schematics, but also an experimental research/teaching platform which is not restricted only to PID based conventional controllers, and can be extended to more advanced control methodologies such as adaptive, robust, real-time, non-linear and/or optimal control.

\subsection{IMU Library}
One library is used to create functions to read from the accelerometer and the gyrometer. It is important to consult the data sheet for each sensor to locate the correct bit address to control sensor. The process begins by sending commands to each sensors unique I2C address: 0x53 for the accelerometer, and 0x69 for the gyrometer. Next, the library writes to different bytes to set the measurement mode and the measurement rates of each sensor. Finally, buffers are created that save accelerometer and gyrometer values in the z, y, and z direction. Compass data also uses the same process but uses a library that was created from Love Electronics.

\begin{figure}[!ht]
 \centering
 \includegraphics[width=3.5truein]{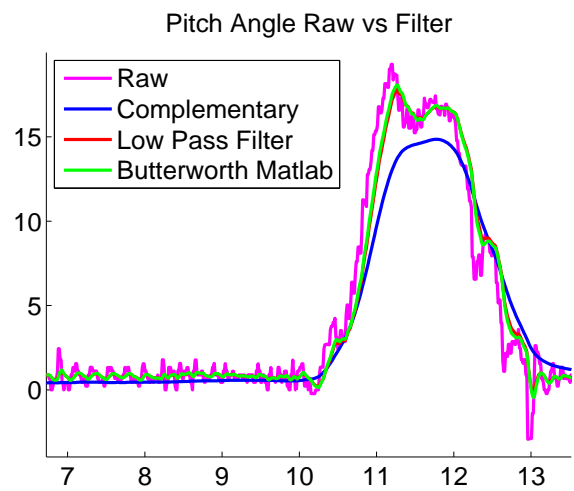}
 \caption{Sample pitch angle filter data} 
 \label{fig:IMUFilter}
 \end{figure}

Initial testing of the values showed very noisy results as shown as a purple line in Figure \ref{Complementary}. Accelerometer readings  have a lot of noise from any small force and a gyrometer is more accurate, but will drift over time when used to calculate the angle states of the quadcopter. Low pass filters and complementary filters are used to remove some of the noise and obtain cleaner measurements for different control schemes. 

The goal of a low pass filter is to allow signals less than the set frequency to pass through which will also remove any high frequency noise. The complementary filter is used to calculate angles (roll and pitch angle of the quadcopter) by combining the best information from the accelerometer and gyrometer. The complementary filter combines the gyroscope data for short term changes combined with the accelerometer for long term data to remove the drift. The general form of the complementary equation is shown in Equation \ref{Complementary}.

\begin{equation}
 angle = \alpha * (angle + gyro * dt) + (1 - \alpha) * accel 
\label{Complementary}
\end{equation}

The complementary filter fuses accelerometer into a first order low pass filter with cut-off frequency of $\alpha$ and integrated gyrometer data into a first order high pass filter of $(1 - \alpha)$ cut-off frequency. Although complementary filtering reduces a large amount of noise, it increases the rise time and the maximum angle never reaches the same value as the raw data shown in purple in Figure \ref{fig:IMUFilter}. Trying improve the performance, low pass filtering (shown in red) and theoretical Butterworth filtering (shown in green) of accelerometer values to calculate the pitch angle were tested. Both filtering techniques showed promising results but the controller stability has yet to be tested. The gyrometer values are filtered using a low pass filter which has a similar form as Equation \ref{Complementary}. The filtered values are computed constantly and used in the PID control scheme.

\subsection{GPS Library}
Data from the GPS is read using the TinyGPS library which is an open source library created by Arduiniana. Most GPS receivers will send out data in standard National Marine Electronics Association (NMEA) sentences at 4800 bytes/s. These NMEA sentences are then parsed by the TinyGPS library to extract information about the time, latitude, longitude, altitude, etc... which will be used in waypoint navigation.

To apply waypoint navigation, first we have to have a destination which will be set by the user in latitude and longitude values. Next, the distance and direction to the waypoint is calculated from the current position. There are two commands in TinyGPS which will be utilized for these calculations: \textit{courseTo()} and \textit{distanceBetween()}. The \textit{courseTo()} command returns the angle from 0 to 360 degrees (using North as a reference of 0 degrees) from the current position to the waypoint position. This information is fused with the compass angle position to find the angle that the quadcopter needs to turn to face the waypoint directly. The \textit{distanceBetween()} command estimates the distance between the current position and waypoint position using the great-circle formula (assumes earth as a hypothetical sphere). This provides enough information to design a control algorithm for waypoint navigation. The GPS currently runs Kalman filter and it is being tested to improve the accuracy of the GPS measurements.

\subsection{Receiver Mapping}\label{PWMmapping}
The receiver values are mapped by triggering interrupts to read in the incomming signal. This is done by setting the 6 analog pins A8-13 to INPUT mode and HIGH to configure the internal pull-up resistors. The same analog pins are set to enable interrupt flags in the PCICR register (PCIE2 for these analog pins) and enabling the individual pins in the PCMSK2 registers. After these setup steps, the corresponding interrupt service routine (ISR(PCINT2\_vect)) will generate an interrupt for every value change from the receiver. 

The time between the rising and the falling edge of the signal will be used to calculate the PWM signal duration. The rising edge is when the signal is first sent out (pin set to HIGH) and the falling edge refers to the end of the signal(pin set to LOW).  In this ISR, the incomming signal is measured from 1024 to 2048 microseconds which corresponds to a PWM signal being sent from the controller input.  One limitation of this is that the Arduino has a resolution of 4 microseconds. This means that the noise in PWM measurement is always a increase/decrease by a factor of 4.

\begin{figure}[!ht]
 \centering
 \includegraphics[width=3.5truein]{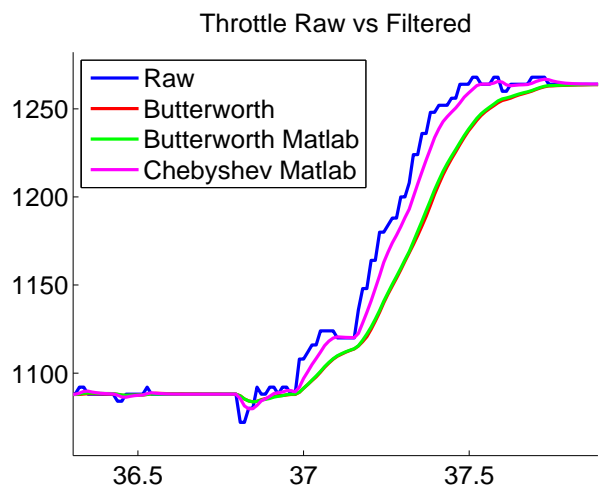}
 \caption{Sample throttle inputs and filter designs} 
 \label{fig:ThrottleFilter}
 \end{figure}
 
Some first order filtering techniques were applied (low pass, Bessel, Chebyshev, etc...) to reduce the noise but it also increased the time delay as seen in Figure \ref{fig:ThrottleFilter}. The Butterworth filter designed in MATLAB (shown in green) matches with the Butterworth filter that is implemented on the quadcopter system (shown in red). Figure \ref{fig:ThrottleFilter} also shows that the Chebyshev filter (shown in purple) has less rise time than the Butterworth but it is at the cost of higher noise.

None of the filtered values are used currently because they would hinder the input reaction speed of a pilot. The receiver signals are mapped from the original (1024 - 2048) values to minimum and maximum motor command values. These values ended up being $\pm 45$ for roll and pitch, and $\pm 135$ for yaw rate (this is 3 times greater than roll and pitch because the pilot desired greater yaw control). Once mapped correctly, they are fed through the control as an input.

\subsection{PID Library}
The PID utilizes the open source PID library created by Beauregard. It is based off the PID equation which can be seen in Equation \ref{equation.PID}. The main advantage of this library is that it already has commands to compute the PID value, set the output limits, set sample time, and change PID gains at any time. Two modifications were made to the library to better fit our needs. 

First, an integral error reset was added to prevent integral windup. The integral term of the PID would continuously increase to the limit before takeoff causing some unnecessary corrections. Second, a negative sign was added to the derivative term of the PID. The library could not compute PID correctly if any of the gains were negative and in this case, the derivative term was negative. This change can be reverted depending on the designed values of the PID controller.

To set up the library in the main file, the PID objects are created for every instance required: roll angle and rate, pitch angle and rate, and yaw rate. This is done by using the constructer function and by providing the values necessary for calculation in the set order. For example, a roll angle PID requires the roll angle measurement, the output roll PID value, the commanded roll angle from the DX6i controller, Kp, Ki, and Kd, values, and finally a direction. The direction will either be DIRECT (positive input to positive output), or REVERSE (negative input to positive output). Once all these constructers are setup correctly, the main script will call on the them to calculate actual input values in the main loop.

\subsection{Chirp Signal Excitation}
A chirp sweep is designed to excite different frequency modes of the quadcopter dynamics. More details about chirp is discussed in section \ref{DataCollection}. Chirp signals were initially done manually, but in the effort to be more efficient, a chirp signal was coded into the main script. The equations used in the code are based off Equations \ref{ChirpEqn} - \ref{ChirpEqn3} which can be found in the book Aircraft and Rotorcraft System Identification \cite{tischler2012}. 

\begin{equation}
 K = C_2 [ e^{ (C_1 * t / T_{rec} )} - 1 ]
\label{ChirpEqn}
\end{equation}

\begin{equation}
 \omega = \omega_{min} + K * ( \omega_{max} - \omega_{min} ) 
\label{ChirpEqn2}
\end{equation}

\begin{equation}
 \delta_{sweep} = A * sin( \omega * t)
\label{ChirpEqn3}
\end{equation}
 
\begin{figure}[!ht]
 \centering
 \includegraphics[width=3.5truein]{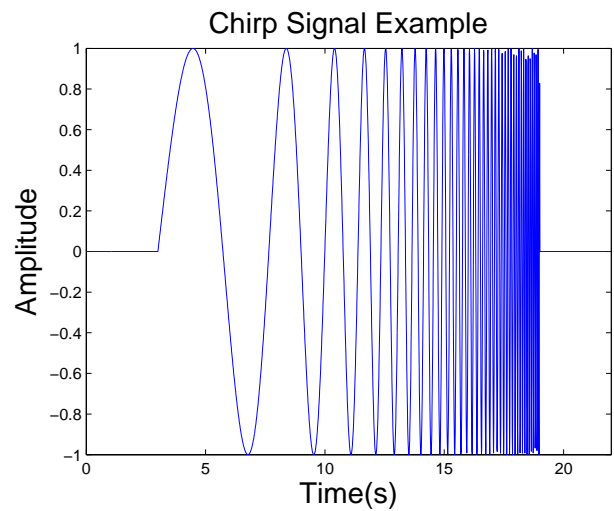}
 \caption{Sample Chirp signal input} 
 \label{fig:ChirpEx}
 \end{figure}

In these equations, $C_1$ and $C_2$ are constants that are defined as 4.0 and 0.0187 respectively from the book. $\omega_{max}$ and $\omega_{min}$ are the minimum and maximum frequencies to sweep through. A is the sweep amplitude constant (which is usually around 10\% of the maximum limit) and t represents time. These equations in combination ared used to create an exponential sweep which spends more time in the lower frequency which is expected to provide better data than using a linear sweep. The equations are coded using Euler integration for $\omega * t$. The chirp signal is set to begin when a switch is changed on the DX6i controller. One important note for running the chirp is that the pilot should try to keep the quadcopter at trim for 3 seconds before the chirp begins and after the chirp ends\cite{tischler2012}. A sample sweep input signal is shown in Figure \ref{fig:ChirpEx}.

\subsection{XBee Data Logging} 
Data logging using the Xbee doesn't require any special library. The Xbee sends and receives data through serial (RX1 and TX1 ports) using serial print statements. Some testing found that data packets would be lost fairly frequently depending on the range of the quadcopter to the laptop.  There were many times when packets of data were lost in flight testing because of signal inteference in the atmosphere that cannot be controlled. Since the loss of data packets through the use of Xbee could not be controlled, all data logging processes now use the micro SD card. The Xbee modules can sbe transitioned into real time telemetry if necesarry in the future.

\subsection{SD Card Data Logging}
The programming for using a micro SD card ended up being a difficult programming task. There are many tutorials that can be found online about how to use the SD library to write data to a SD card, but they require high processing time due to opening and closing the file every loop. Using the SD library was a huge problem since it lowered the sample rate to 20[Hz] with the goal of 100[Hz]. Running the system at 20[Hz] caused the quadcopter to become unstable. A solution sample rate issue was found using the basic SdFat library. The SD library is a wrapper of the SdFat library, which means higher performance can be achieve at the expense of coding complexity using the SdFat library. 

The code first creates a contiguous file in the setup loop using the SdFat command \textit{createContiguous()}. This will set up a file that all the data is saved to the same place for quicker access. The naming of this file will only accept 12 characters total. This leaves a total of 8 characters to name a file because the file extension will take up 4 characters (.txt, .csv, etc.. counts as 4 characters). A file is setup with multiple 512 byte blocks then collects data in the main loop as 512 byte data strings. This requires converting some of the double and float values into string format using the \textit{dtostrf()} function. The string is then transfered to a uint8\_t buffer, which stands for a unsigned integer of length of 8 bits (this converts the string to the byte format). Once in a uint8\_t format, the buffer is written to a 512 byte block in the file that was set up at the start. This process was found to save enough time to keep the sample rate over 100[Hz].

\subsection{Simulink Autocoding}
MATLAB Simulink(c) is a widely used program that has capabilities to auto generate Arduino code from a model. This allows a user the ability to design a controller in the Simulink environment, generate autocode, and upload it to the Arduino more easily than translating the algorithm into Arduino code manually. To achieve this capability, the Simulink S-functions are used to allow creation of Simulink blocks which mimic Arduino code. These S-functions act as user defined functions which allow each module to interact with the Simulink environment. 

The process begins by creating S-function blocks for each module of the quadcopter: receiver, IMU, ESC and motors, micro SD card, etc.. This is done by placing the current Arduino code  into the correct section of the S-function setup. The Data Properties tab holds all named  variables that are to be received or sent (input and output signals like a normal Simulink block). The Libraries tab holds all the associated Arduino libraries, setup variables declarations, and any object creation instances. The Discrete Update contains setup functions that are usually in the Arduino void setup loop (anything that needs to only be run one time. The Outputs tab contains the main code that is looped (similar to Arduino void loop). Once all the S-functions are complete, they are compiled together and the model can be built directly to the Arduino Mega.

\begin{figure*}[!ht]
 \centering
 \includegraphics[width=7.0truein]{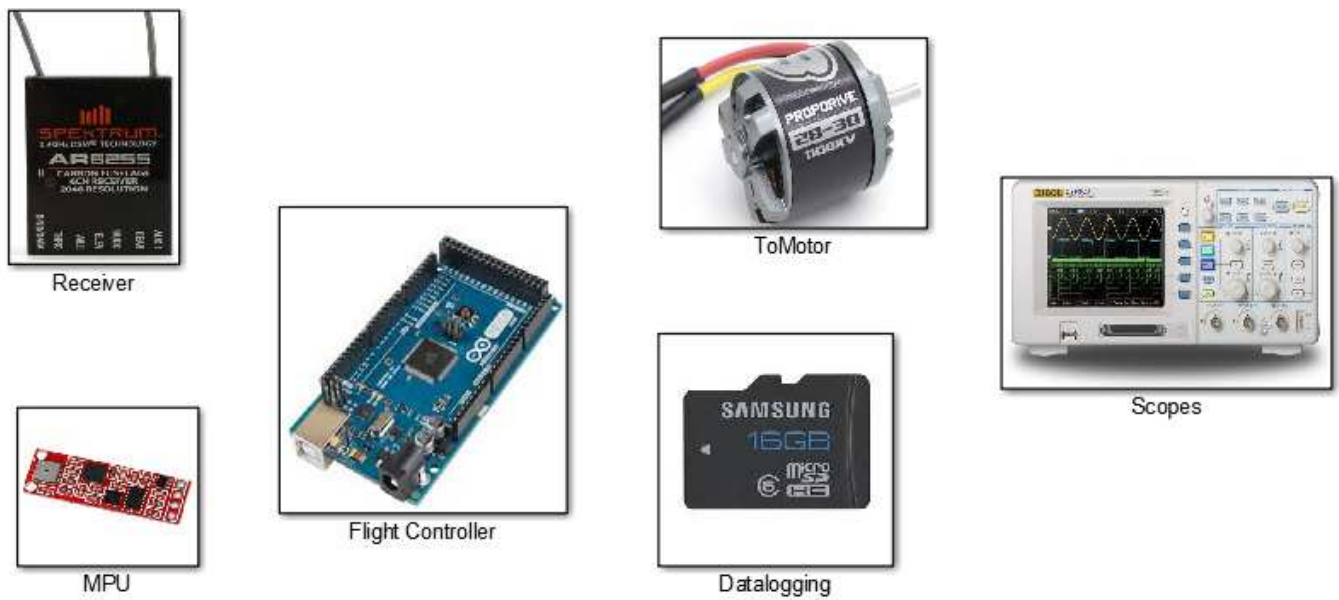} 
 \caption{Quadcopter hardware in the loop Simulink file} 
 \label{fig:Autocoding}
 \end{figure*}
 
\begin{figure*}[!ht]
 \centering
 \includegraphics[width=7.0truein]{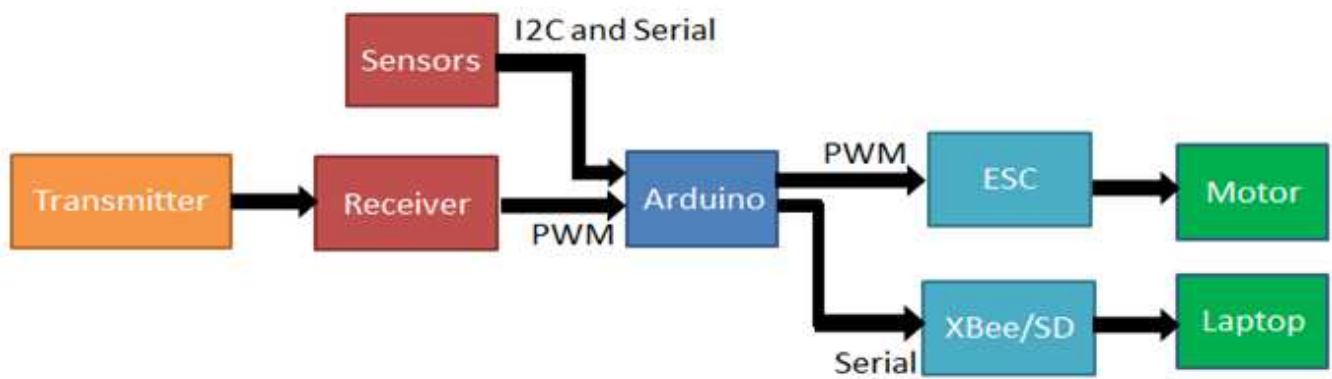}
 \caption{Flow Chart for Arduino Modules} 
 \label{fig:FlowChart}
 \end{figure*} 

Figure \ref{fig:Autocoding} shows what the current quadcopter Simulink file which is modeled in a flowchart shown in Figure \ref{fig:FlowChart}. Each block contains an S-function that performs the same task it does on the Arduino Mega. One exception to this is the Scopes block which contains scopes to be viewed in real time external mode. The other exception is the Flight Controller block which contains PID blocks instead of using the Arduino PID library. This PID block is used to create the same functionality as the Arduino code using Simulink blocks to verify that both cases demonstrate the same behavior. This file can be loaded onto the Arduino but there have been some issues with sampling time and stability which require further testing and debugging.

\section{State Space Model and Stability Analysis} \label{SSMSA}
A state space matrix representation of the quadcopter can be derived using equations of motion discussed in Section \ref{sec:Modeling}. This model will be used to calculate initial PID values to stabilize the quadcopter platform. Since the model is a rough estimate these PID values, which only are starting estimates, are further tuned to accomodate pilot's preference and mission's requirements.

\subsection{State Space model} \label{SSM}
The linearized state space model (SSM) is derived using the Simulink block model discussed in Section \ref{sec:Modeling}, where the input values for the quadcopter were provided in Table \ref{table.parameters}. The final state space matrices are shown in Eq.'s \ref{Eqn.AB}-\ref{Eqn.CD}.

\setcounter{MaxMatrixCols}{12}
\begin{equation}\label{Eqn.AB}
A = \begin{bmatrix}
 0 & 0 & 0 & 0 & 0 & 0 & 0 & 0 & 0 & 0 & 9.81 & 0 \\
 0 & 0 & 0 & 0 & 0 & 0 & 0 & 0 & 0 & -9.81 & 0 & 0 \\
 0 & 0 & 0 & 0 & 0 & 0 & 0 & 0 & 0 & -9.81 & 0 & 0 \\
 0 & 0 & 0 & 0 & -1.63 & 0 & 0 & 0 & 0 & 0 & 0 & 0 \\
 0 & 0 & 0 & 1.63 & 0 & 0 & 0 & 0 & 0 & 0 & 0 & 0 \\
 0 & 0 & 0 & 0 & 0 & 0 & 0 & 0 & 0 & 0 & 0 & 0 \\
 1 & 0 & 0 & 0 & 0 & 0 & 0 & 0 & 0 & 0 & 0 & 0 \\
 0 & 1 & 0 & 0 & 0 & 0 & 0 & 0 & 0 & 0 & 0 & 0 \\
 0 & 0 & 1 & 0 & 0 & 0 & 0 & 0 & 0 & 0 & 0 & 0 \\
 0 & 0 & 0 & 1 & 0 & 0 & 0 & 0 & 0 & 0 & 0 & 0 \\
 0 & 0 & 0 & 0 & 1 & 0 & 0 & 0 & 0 & 0 & 0 & 0 \\
 0 & 0 & 0 & 0 & 0 & 1 & 0 & 0 & 0 & 0 & 0 & 0
\end{bmatrix}
\end{equation}
\begin{equation}
B = \begin{bmatrix}
 0 & 0 & 0 & 0 \\
 0 & 0 & 0 & 0 \\
 0.7143 & 0 & 0 & 0 \\
 0 & 12.3457 & 0 & 0 \\
 0 & 0 & 12.3457 & 0 \\
 0 & 0 & 0 & 7.0423 \\
 0 & 0 & 0 & 0 \\
 0 & 0 & 0 & 0 \\
 0 & 0 & 0 & 0 \\
 0 & 0 & 0 & 0 \\
 0 & 0 & 0 & 0 \\
 0 & 0 & 0 & 0
\end{bmatrix},~~~
x = \begin{bmatrix}
     X \\ Y \\ Z \\ U \\ V \\ W \\ \phi \\ \theta \\ \psi \\ P \\ Q \\ R
\end{bmatrix}
\end{equation}

\begin{equation}
C = I_{(12x12)}, 
D=\begin{bmatrix}
0
\end{bmatrix}
\label{Eqn.CD}
\end{equation}

\subsection{Open-loop Analysis}
Analyzing the plant dynamics indicates an initially unstable system, as expected. Thus, any small disturbance will cause the quadcopter to become unstable without a controller to stabilize it. However, the SSM is fully controllable due to the different combination of propeller speed inputs.
Analysis of the state space model shows that disturbances in roll and pitch angles result in very large displacements. Roll and pitch angles increase over 20 degrees in 10 seconds because the roll and pitch rates do not stabilize after an impulse. The yaw angle and rate stay small in comparison to pitch and roll angles. The translational modes of the SSM show some very odd behavior because they exclude aerodynamic forces. The result of the SSM analysis show that the overall system is unstable and translation modes cannot be controlled directly, which brings the need for closed-loop system analysis.

\subsection{Closed-loop Analysis}
An initial Proportional Integral Derivative (PID) inner-loop controller is designed to stabilize the system where the PID gain value characteristics inherits the conventional form as shown in Eq.(\ref{equation.PID}) \cite{sa2012,dicesare2009} .

\begin{equation}
\label{equation.PID}
PID(s) = K_p + \frac{K_i}{s} + T_d s
\end{equation}

\begin{figure*}[!ht]
 \centering
 \includegraphics[width=6.0truein]{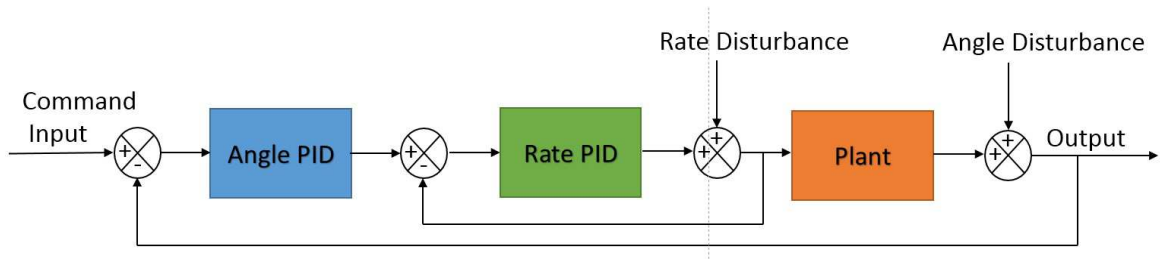}
 \caption{Quadcopter inner and outer loop PID diagram}
 \label{fig:PIDFlowChart}
 \end{figure*}

For stability, a cascaded control architecture will be implemented which uses multiple inner loops and multiple input signals. The innermost PID will use rate values and the outer loop will stabilize any angular disturbances (refer to Figure \ref{fig:PIDFlowChart} ). The PID values chosen to stabilize the open loop dynamics are shown in Table \ref{table.PIDgains}.

\begin{table*}[!ht]
\centering
\caption{Initial PID values to stabilize roll and pitch.}
\begin{tabular}{|l|l|l|l|}
  \hline
  & Roll/Pitch Angle Gains & Roll/Pitch Rate Gains & Yaw Rate Gains \\ \hline
  $K_p$ & 3.604  & 0.2209 & 0.1141\\
  $K_i$ & 0      & 0 	  & 0.6340\\
  $T_d$ & 0      & 0.014  & 0\\
  \hline
\end{tabular}
\label{table.PIDgains}
\end{table*}

\begin{figure*}[htbp!]
 \centering
 \includegraphics[width=6.0truein]{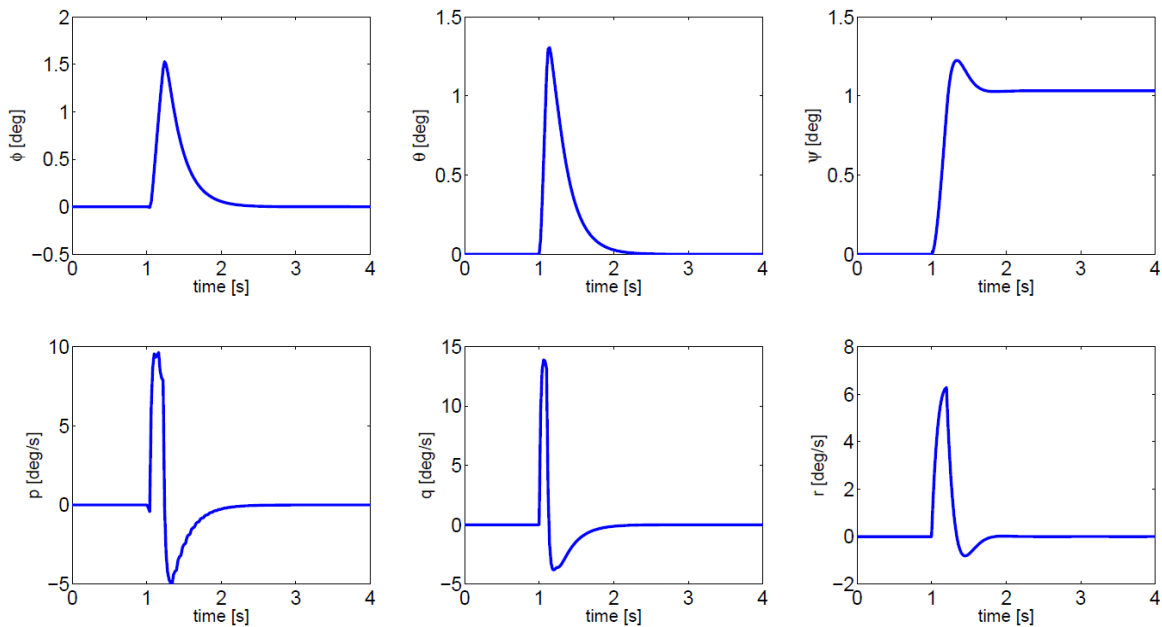}
 \caption{Impulse response of closed loop PID quadcopter dynamics.}
 \label{fig:PID}
 \end{figure*}

The goal is to use these values as a starting point to achieve a stabilized closed-loop system. The PID gains are further tuned to the pilot's preference to perform sweeping maneuvers for system identification. Although the quadcopters EoMs are coupled, PID controllers for the other states (such as translation) are not necessary because the disturbances are relatively small compared to the roll and pitch angles. The closed loop response of the quadcopter system with the PID values from Table \ref{table.PIDgains} are shown in Figure \ref{fig:PID}.

It is clear from Figure \ref{fig:PID} that roll, pitch, and yaw responses are stable and return back to the trim value at zero within a second. Since there is no closed loop control on position, translational motion(x,y,z) of the quadcopter will have small steady state errors, which are relatively small and therefore negligable.  

\section{System Identification} \label{sec:SystemID}
As stated in \cite{miller2007}, \cite{lee2011} and \cite{tischler2012}, it is well known that system identification can be used to develop a linearized SSM using experimental flight data. This involves mapping a known input signal to the flight data response in the frequency domain to obtain the transfer function of the corresponding state. Once this is achieved for all 6 states of the system, the model structure can be built from the estimated transfer functions for further theoretical/simulational analysis. The goal is to compare and verify the system identification model with actual quadcopter experimental data.
 
\subsection{Data Collection} \label{DataCollection}
\begin{figure}[ht!]
 \centering
 \includegraphics[width=3.5truein]{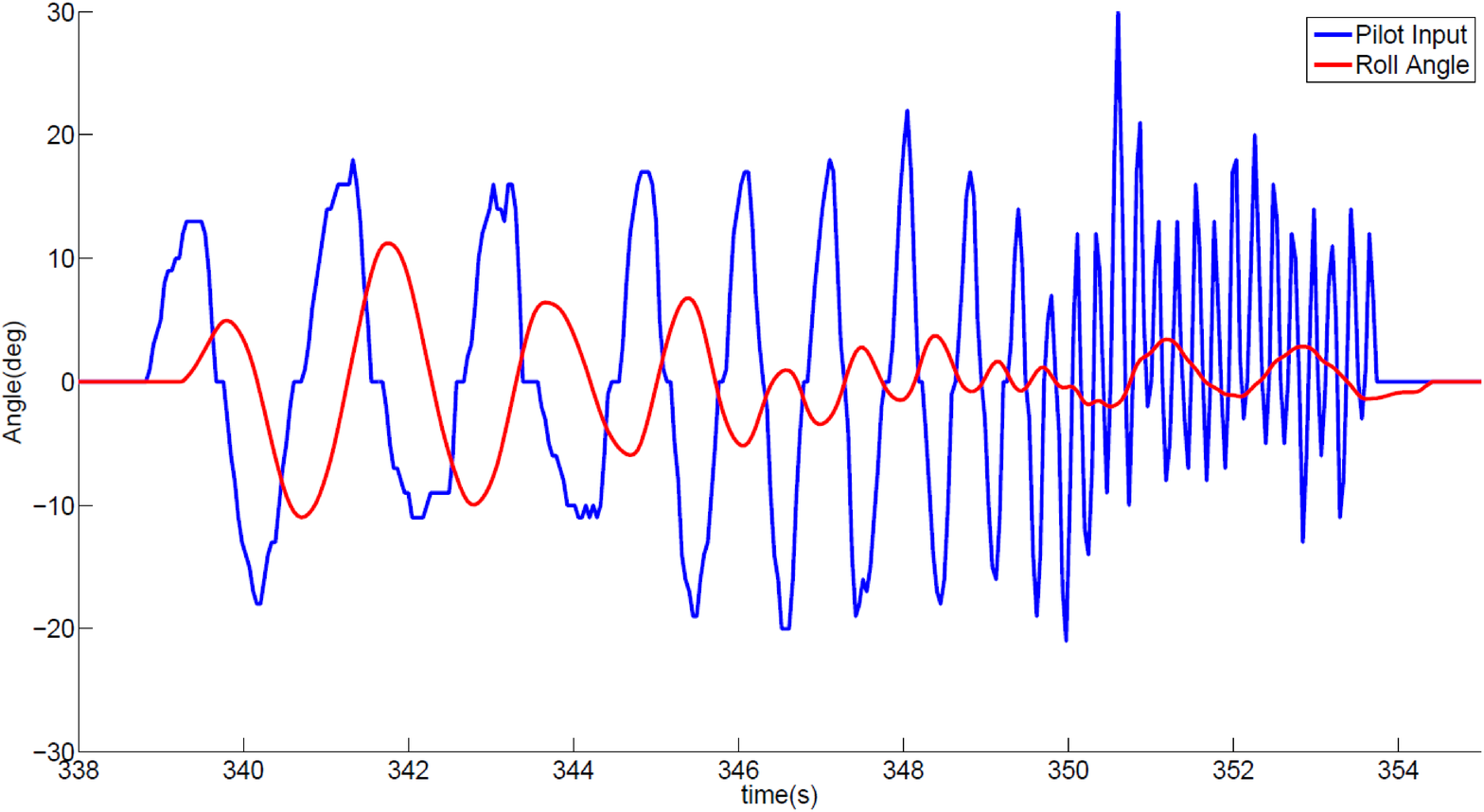}
 \caption{Roll Coherence in $CIFER^{\circledR}$ of a set of sweep data.}
 \label{fig:Sweep}
 \end{figure}

The identification process starts with retrieving frequency domain sweep data via a chirp signal. A chirp signal is essentially a sinusodial function which starts at a low frequency and slowly increases to higher frequencies to cover (and excite) all the different modes of the system. The sweep, in our case, is implemented manually through a pilot input and is applied to roll, pitch, and yaw states with the output data for angles and rates being logged for analysis in \cite{cifer2006}. A sample sweep in roll from the pilot input is shown in Figure \ref{fig:Sweep}. The magnitude and frequency of the sweep command vary since a pilot input is not perfect, but it still provides valuable data for analysis.

\subsection{SISO Identification Analysis}
For system identification analysis, the student version of \cite{cifer2006} (Comprehensive Identification from Frequency Response) is utilized . It is an identification software that is designed to use time domain test data to extract frequency domain models in the form of transfer functions and/or state space models. In this study, only angular modes were able to be analyzed because the translational measurements were not precise enough for analysis. The resulting identified SISO transfer functions for roll and pitch are shown in Eq.'s (\ref{equation.RollTF}) - (\ref{equation.PitchTF}). 

\begin{equation}
\begin{split}
\label{equation.RollTF}
TF_{RollAngle}=& \frac{2.305s}{s^2 + 3.894s + 3.967} * e^{-\tau_1 s} 
\end{split}
\end{equation}

\begin{equation}
\begin{split}
\label{equation.PitchTF}
TF_{PitchAngle}=& \frac{2.008s}{s^2 + 3.206s + 4.058} * e^{-\tau_2s}
\end{split}
\end{equation}

\begin{equation}
\begin{split}
\label{equation.YawTF}
TF_{YawRate}=& \frac{10.68s + 138.8}{s^2 + 11.64s + 163.8} * e^{-\tau_3s}
\end{split}
\end{equation}
where associated time delays are
\begin{equation}
\tau_1 = 0.197[sec],~~ \tau_2 = 0.2[sec],~~ \tau_3 = 0.0592[sec]
\end{equation}

\begin{figure*}[!ht]
 \centering
 \includegraphics[width=7truein]{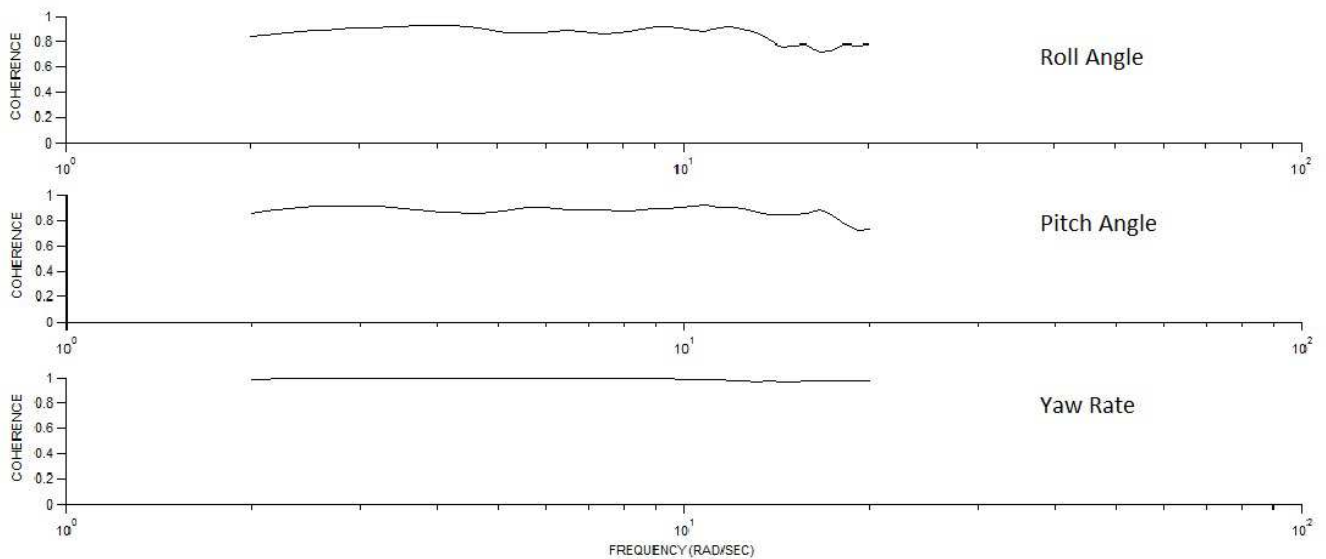}
 \caption{Coherence from $CIFER^{\circledR}$ using automated sweep data.}
 \label{fig:Coherence}
 \end{figure*}

\begin{figure*}[!ht]
 \centering
 \includegraphics[width=7truein]{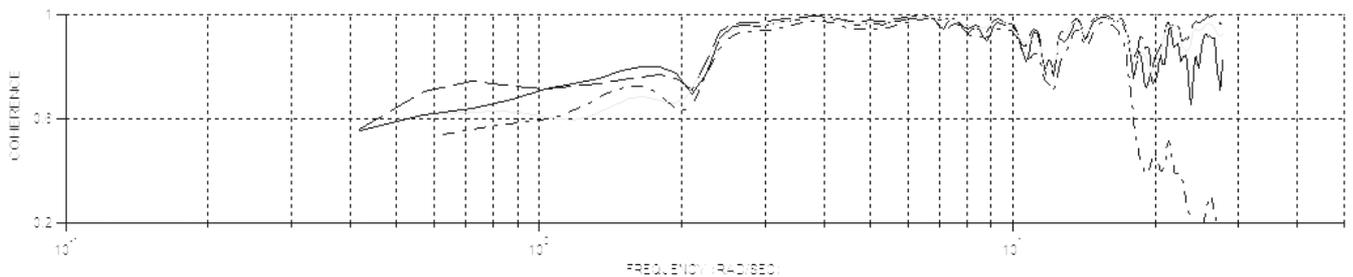}
 \caption{Roll angle coherence from $CIFER^{\circledR}$ using manual sweep data.}
 \label{fig:CoherenceManual}
 \end{figure*}

The transfer functions obtained from automated sweeps are acceptable as determined by the coherence shown in Figure \ref{fig:Coherence}. All three figures show that the coherence varies around 0.8 to 1. This also shows a large improvement in using automated chirps compared to using manual pilot chirps. The manual chirp coherence for roll can be seen in Figure \ref{fig:CoherenceManual} which shows that lower frequencies are adequate but at higher frequencies, the sweep data is not consistent enough past $10 [rad/s]$. The coefficient in the exponential term represents the time delay from the input to output.

\section{Time Domain Validation}\label{sec:validation}
\begin{figure}[!ht]
 \centering
 \includegraphics[width=3.5truein]{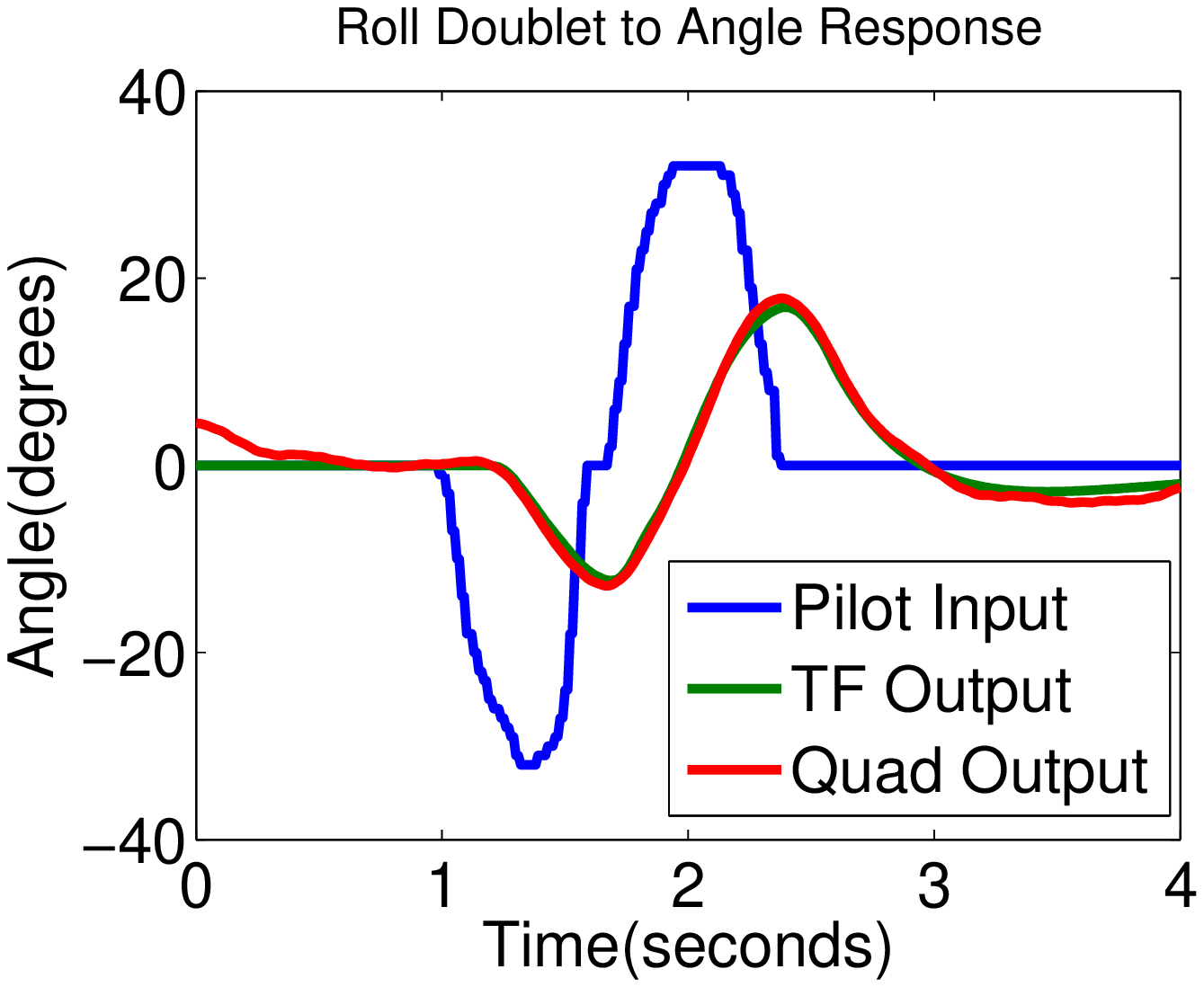}
 \includegraphics[width=3.5truein]{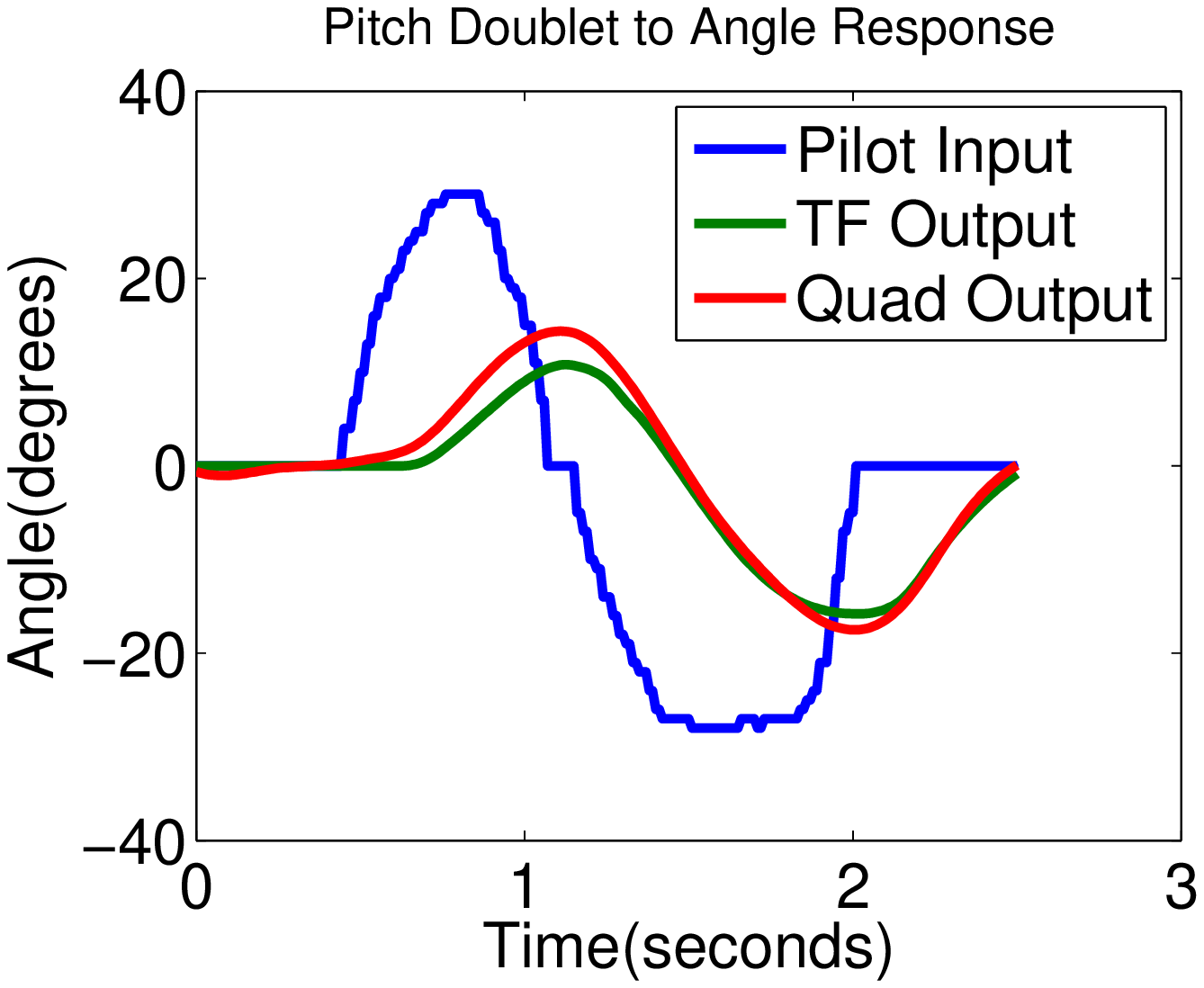}  \includegraphics[width=3.5truein]{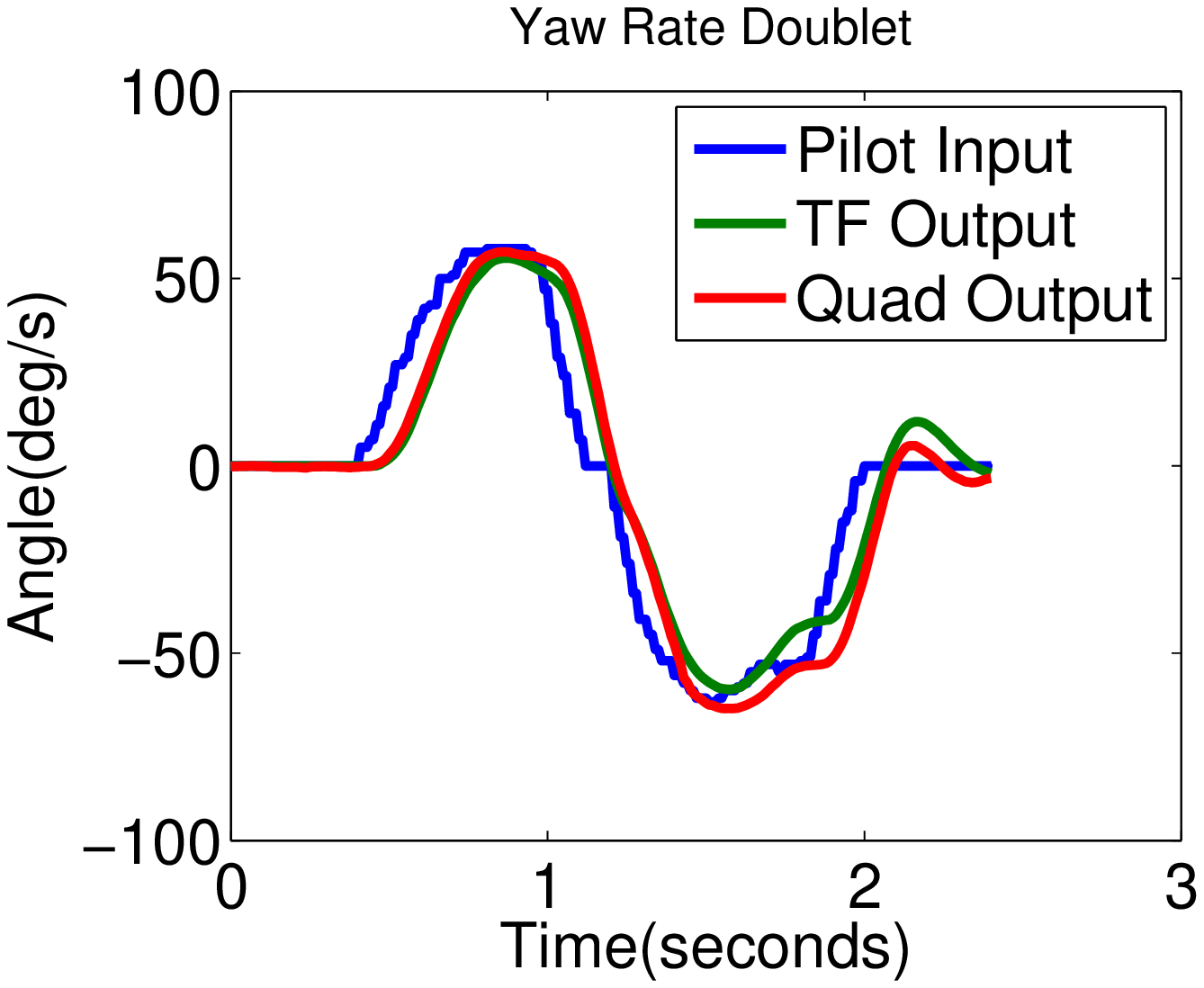}
 \caption{Doublet response of quadcopter compared with $CIFER^{\circledR}$ transfer function.}
 \label{fig:Doublet}
 \end{figure}

Following to identification results, derived dynamics are compared with the actual flight data to verify the accuracy of the model dynamics. The input signals are manual piloted doublets and the results are compared with actual flight data. The result of roll and pitch doublets are shown in Figure \ref{fig:Doublet} using the transfer function from Equations \ref{equation.RollTF} and \ref{equation.PitchTF}.

It can be seen that the time delay of the transfer function matches up with the actual system very well. The roll angle transfer function doublet magnitude also matches up with flight data but the pitch angle transfer function does reach the correct magnitude, even though the high coherence implies that they should match. The reason for lower magnitude is because the transfer function that $CIFER^{\circledR}$ extracts is a lower-order equivalent system (LOES). This SISO transfer function can be used as an initial model, but it does not account for the quadcopter dynamics that may be of a higher order and involve coupled dynamic behavior.
 
A complete SSM is desired so that it can be used as for analytic modeling and verification. However, a complete SSM could not be extracted due to a lack of consistent measurements. The problem lies in gathering correlatable off-axis chirp data. Identification of a SSM model requires accurate measurements for off-axis state to account for the coupling dynamics. It can be seen from the previous section that on-axis data has high coherence, but the off-axis coherence is contaminated with noise, or woud have low coherence at certain frequencies due to using filtered data. The only set of data that $CIFER^{\circledR}$ did not error on was the yaw input to yaw rate output which had low coherence in the 20s. Different filter values and chirp magnitudes were tested for multiple flight tests but almost all trials resulted in a low coherence in the SSM identification.

\section{Conclusion}\label{sec:conclusion}
In this study, we provided a process of assembling an autonomous, low-cost, off-the-shelf product based quadcopter platform from scratch. We designed control laws to stabilize it using an Arduino Mega. The EoMs for a quadcopter have been studied and a linear state space model was extracted for analysis. The quadcopter testbed was assembled and programmed using the Arduino Mega to be flown using a DX-6 controller. One novelty we presented in this paper is that we provide a modified, novel version of Arduino code which is not restricted to PID controllers only, and can be extended to more advanced control methodologies due to its modular and highly flexible structure. For state measurements, the IMU and GPS have been calibrated with a many different types of filters. A closed loop PID is designed in Simulink using the state space model to stabilize non-linear plant dynamics. The closed loop controller has also been coded into the Arduino Mega. With this design, flight data is recorded on the micro SD card for analysis. Transfer functions for roll and pitch were identified and verified using doublet data from real flight. SSM extraction becomes problematic due to problems with low signal-to-noise ratio and low coherence of off-axis data, and is currently and on-going search where the results will be reported in another study.


\end{document}